\documentclass[journal]{AIAA}

\usepackage[utf8]{inputenc}
\usepackage{graphicx,float,placeins,tikz,onimage,caption,subcaption,xcolor}
\usepackage{physics,amsmath,mathtools,amsfonts,mathrsfs,bm}

\let\oldopenbox\openbox
\let\openbox\relax
\usepackage{amsthm}
\let\openbox\oldopenbox

\let\oldBbbk\Bbbk
\let\Bbbk\relax
\usepackage{amssymb}
\let\Bbbk\oldBbbk

\usepackage{longtable,tabularx,booktabs,multirow}
\usepackage[version=4]{mhchem}
\usepackage{siunitx}
\usepackage{textcomp}
\usepackage[pagewise]{lineno}
\usepackage{wasysym}
\usepackage{ifthen}
\usepackage{enumitem}
\usepackage{url}

\PassOptionsToPackage{colorinlistoftodos,prependcaption,textsize=normalsize}{todonotes}
\usepackage{todonotes}

\usepackage{listings}
\lstdefinestyle{codestyle}{
    commentstyle=\color{codegreen},
    keywordstyle=\color{blue},
    numberstyle=\tiny\color{gray},
    stringstyle=\color{codepurple},
    basicstyle=\fontfamily{pcr}\selectfont\footnotesize,
    breakatwhitespace=false,
    breaklines=true,
    captionpos=b,
    keepspaces=true,
    numbers=left,
    numbersep=5pt,
    showspaces=false,
    showstringspaces=false,
    showtabs=false,
    tabsize=2
}
\newcommand{\reals}{\mathbb{R}}

\newcommand{\integers}{\mathbb{Z}}

\newcommand{\mF}{\mathcal{F}}
\renewcommand{\P}[1][P]{\mathbb{#1}}

\definecolor{White}{RGB}{255, 255, 255}
\definecolor{Silver}{RGB}{192, 192, 192}
\definecolor{Grey}{RGB}{128, 128, 128}
\definecolor{Black}{RGB}{0, 0, 0}
\definecolor{Tomato}{RGB}{255, 99, 71}
\definecolor{Red}{RGB}{255, 0, 0}
\definecolor{Maroon}{RGB}{128, 0, 0}
\definecolor{Magenta}{RGB}{172, 87, 156}
\definecolor{Yellow}{RGB}{255, 255, 0}
\definecolor{Tan}{RGB}{205, 181, 145}
\definecolor{Olive}{RGB}{128, 128, 0}
\definecolor{Lime}{RGB}{0, 255, 0}
\definecolor{Green}{RGB}{0, 128, 0}
\definecolor{Chartreuse}{RGB}{160, 252, 78}
\definecolor{Aqua}{RGB}{0, 255, 255}
\definecolor{SkyBlue}{RGB}{83, 194, 236}
\definecolor{Teal}{RGB}{0, 128, 128}
\definecolor{Blue}{RGB}{0, 0, 255}
\definecolor{DodgerBlue}{RGB}{85, 154, 248}
\definecolor{RoyalBlue}{RGB}{65, 105, 225}
\definecolor{NavyBlue}{RGB}{25, 25, 107}
\definecolor{TealBlue}{RGB}{84, 113, 171}
\definecolor{Navy}{RGB}{0, 0, 128}
\definecolor{Fuchsia}{RGB}{232, 51, 244}
\definecolor{Purple}{RGB}{128, 0, 128}
\definecolor{Orange}{RGB}{242, 132, 52}

\definecolor{colunstablethreefour}{RGB}{248, 130, 83}
\definecolor{colstablethreefour}{RGB}{136, 183, 213}
\definecolor{colunstablefivesix}{RGB}{166, 135, 186}
\definecolor{colstablefivesix}{RGB}{130, 186, 126}

\definecolor{codegreen}{RGB}{0, 153, 0}
\definecolor{codepurple}{RGB}{148, 0, 209}

\definecolor{python_deep_blue}{RGB}{27, 106, 148}
\definecolor{python_deep_orange}{RGB}{228, 116, 50}
\newboolean{includefigures}
\newboolean{revisionmarkups}

\usepackage[normalem]{ulem}
\useunder{\uline}{\ul}{}

\graphicspath{{"Figures"}}
\setboolean{includefigures}{true}
\setboolean{revisionmarkups}{false}

\ifthenelse{\boolean{revisionmarkups}}
{
    \usepackage[authormarkup=none, commentmarkup=todo, todonotes={textsize=tiny}]{changes}
    \definechangesauthor[name=Amlan Sinha]{AS}
}{
    \usepackage[final]{changes}
}

\newcommand{\orcid}[1]{\href{https://orcid.org/#1}{\includegraphics[scale=0.25]{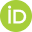}}}

\title{Initial Guess Generation for Low-Thrust Trajectory Design with Robustness to Missed-Thrust-Events}

\author{
Amlan Sinha 
\orcid{0009-0001-4572-7583} 
\footnote{Corresponding Author: \href{mailto:amlans@princeton.edu}{amlans@princeton.edu}} 
\footnote{
Ph.D. Candidate, Mechanical and Aerospace Engineering, Princeton University, Princeton, NJ 08540.
Member AIAA.
} 
and 
Ryne Beeson 
\orcid{0000-0003-2176-0976} 
\footnote{
Assistant Professor, Department of Mechanical and Aerospace Engineering, Princeton University, Princeton, NJ 08540.
Member AIAA.
}
}

\affil{Princeton University, Princeton, NJ 08840}

\begin{document}

\maketitle


\begin{abstract} 
The growing interest in cislunar space exploration in recent years has driven an increasing demand for efficient low-thrust missions to key cislunar orbits. 
These missions, typically possessing long thrust arcs, are particularly susceptible to operational uncertainties such as missed thrust events. 
Addressing these challenges requires efficient robust trajectory design frameworks during the preliminary mission design phase, where it is necessary to explore the solution space at a rapid cadence under evolving operational constraints.
However, existing methods for missed thrust design rely on solving high-dimensional nonlinear programs, where generating effective initial guesses becomes challenging.
To enhance computational efficiency, quality, and depth of robustness of solutions from global search, we compare two initial guess strategies: a baseline \emph{non-conditional} global search, which samples from a static distribution with global support, and a \emph{conditional} global search, which generates initial guesses conditioned on solutions to problems with less depth of robustness. 
The conditional search provides a sequential procedure for solving increasingly robust problems. 
We validate the improvements in the conditional approach using a low-thrust case study for the Lunar Gateway Power and Propulsion Element, where our results demonstrate that it significantly improves convergence rate and solution quality, highlighting its potential in preliminary robust trajectory design.
\end{abstract}
\section*{Nomenclature}

{\renewcommand\arraystretch{1.0}
\noindent\begin{longtable*}{@{}l @{\quad=\quad} l@{}}

$\Omega$                & random sample space \\
$\mF$                   & $\sigma$-algebra of measurable events \\
$(\mF_t)_{t \geq 0}$    & filtration representing information over time \\
$\P$                    & probability measure \\
$\omega$                & random variable \\
$J$                     & objective function \\
$\phi$                  & terminal cost \\
$\mathcal{L}$           & running cost \\
$\bm{f}$                & uncontrolled state vector field \\
$\bm{g}$                & controlled state vector field \\
$t_0$                   & initial time \\
$t_f$                   & final time \\
$\bm{\xi}$              & spacecraft state \\
$\bm{u}$                & spacecraft control \\ 
$\bm{x}$                & nonlinear programming decision variable \\
$\bm{c}$                & nonlinear programming constraint \\
$T_s$                   & shooting time \\
$T_i$                   & initial coast time \\
$T_f$                   & final coast time \\
$\mathcal{E}$           & equality index set \\
$\mathcal{I}$           & inequality index set \\
$N$                     & number of segments \\ 
$\overline{N}$          & number of decision variables \\ 
$\tau_1$                & time where a missed thrust event begins \\
$\Delta \tau$           & duration of the missed thrust event \\
$\mathcal{X}$           & solution space \\
$X$                     & initial guess generator \\
$\mathcal{M}$           & mapping strategy \\
$\mathcal{S}(k)$        & Non-conditional strategy for a robust problem with $k$ realizations \\
$\mathcal{S}(k | k')$   & Conditional strategy for a $k$-robust problem with initial guesses derived from a $k'$-robust problem.\\
\multicolumn{2}{@{}l}{Superscripts}\\
$\dagger$       & reference solution \\
$\omega$        & realization solution \\
\multicolumn{2}{@{}l}{Subscripts}\\
$i, j$          & indices \\
$k, k'$         & number/index of realizations \\
\end{longtable*}}
\section{Introduction}
\label{section: introduction}

\lettrine{I}{n} the recent past, low-thrust (LT) propulsion systems are becoming increasingly popular for space missions, particularly cislunar. 
Several Artemis I payloads, such as Lunar IceCube, EQUULEUS, Near-Earth Asteroid Scout, BioSentinel, LunIR, and Lunar Polar Hydrogen Mapper, exemplify this trend.
Despite their many advantages, LT spacecrafts are particularly susceptible to \emph{safe mode events}, which can occur if an anomalous event (e.g., impact with space debris) forces the spacecraft to enter a protective mode during which all thruster operations are switched off. 
If such an event occur coincides with a scheduled thrust arc, it results in what is known as a \emph{missed thrust event} (MTE). 
Due to their characteristically long thrust arcs, MTEs are relatively common in LT missions \cite{imken_modeling_2018}.
If not adequately addressed during the preliminary design phase, they can severely compromise mission performance and, in certain cases, lead to complete failure - especially if they disrupt maneuvers which must occur at critical junctures along the trajectory (e.g., flybys).

\ifthenelse{\boolean{includefigures}}
{
    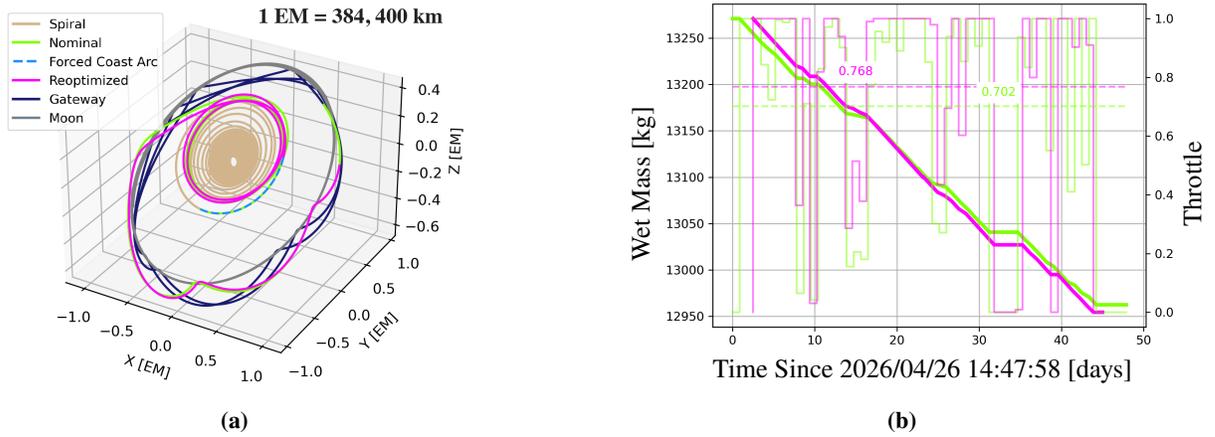
\begin{figure}[!htb]
        \centering
        \begin{subfigure}[b]{0.475\textwidth}
            \centering
            \includegraphics[keepaspectratio, width=0.8\textwidth]{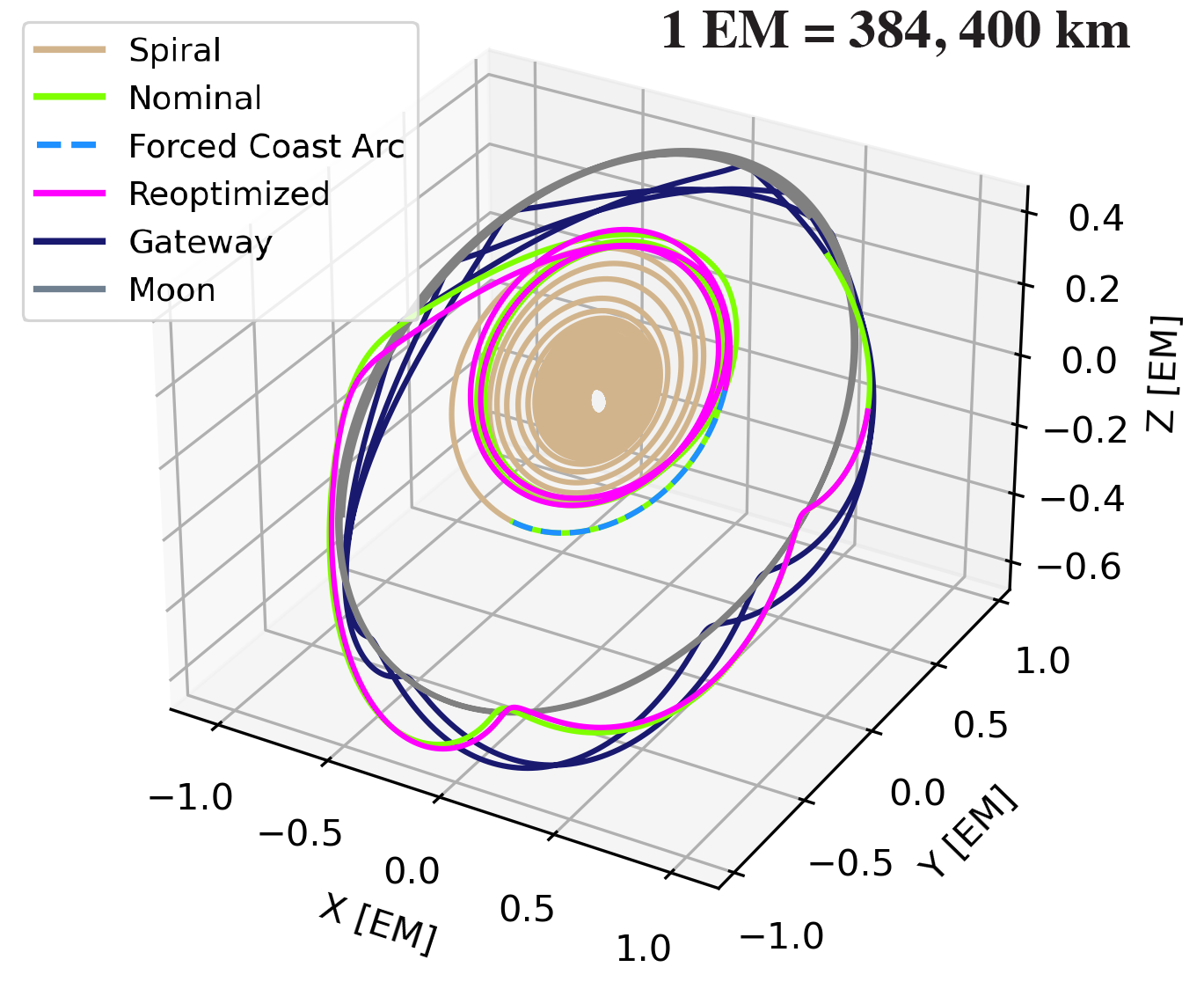}
            \caption{}
            \label{fig: motivating_trajectory}
        \end{subfigure}
        \hfill
        \begin{subfigure}[b]{0.45\textwidth}
        \centering
        \begin{tikzpicture}
            \node[anchor=south west,inner sep=0] (image) at (0,0) {\includegraphics[keepaspectratio, width=\textwidth]{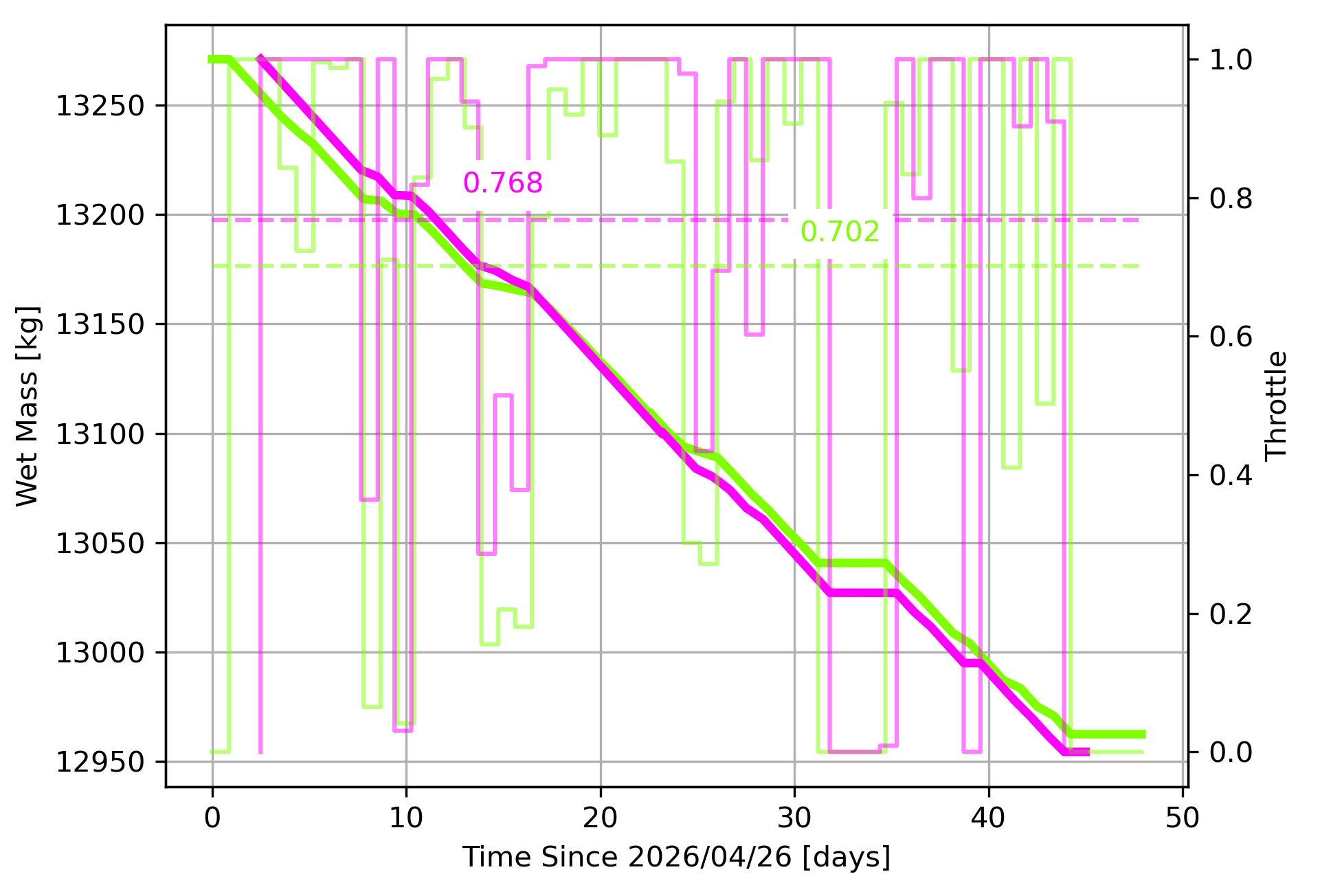}};
            \begin{scope}[x={(image.south east)},y={(image.north west)}]
                \node[fill=white, opacity=1.0, text opacity=1, anchor=south] at (0.5,-0.05) {Time Since 2026/04/26 14:47:58 [days]};
                \node[fill=white, opacity=1.0, text opacity=1, anchor=south, rotate=90] at (0.04,0.5) {Wet Mass [kg]};
                \node[fill=white, opacity=1.0, text opacity=1, anchor=south, rotate=90] at (1.015,0.5) {Throttle};
            \end{scope}
        \end{tikzpicture}
        \caption{}
        \label{fig: motivating_mass_throttle}
    \end{subfigure}
        \caption{
        Reoptimizing a nominal solution following an MTE can deteriorate the baseline
        }
        \label{fig: motivating}
    \end{figure}
}
{
}

As a motivating example, consider a LT spacecraft, modeled after the Power and Propulsion Element \cite{mcguire_power_2020, mcguire_overview_2021} (discussed in more detail in \S \ref{section: problem setup}), on a spiral trajectory in the cislunar realm, where the objective is to identify a feasible, minimum fuel transfer to the Earth–Moon $\mathcal{L}_2$ Southern 9:2 Near Rectilinear Halo Orbit.
Due to the low thrust acceleration, solutions to this problem typically span multiple revolutions with long thrust arcs, making them particularly vulnerable to MTEs.
For illustration, let us consider a nominal solution which consumes approximately 2,038 kg of propellant over a total flight time of 225 days (Fig. \ref{fig: motivating}). 
To understand how an MTE could affect this solution, we simulate a 60-hour outage near the beginning of the transfer, preventing the spacecraft from executing the first scheduled burn, and forcing it to coast under natural dynamics until it can resume its thruster operations once again. 
Reoptimizing this solution after the recovery yields a trajectory which differs significantly from the original.

Figure \ref{fig: motivating_trajectory} illustrates the nominal non-robust trajectory in \textcolor{Chartreuse}{\textbf{green}} alongside the reoptimized trajectory in \textcolor{Fuchsia}{\textbf{pink}}. 
The LT spiral phase is highlighted in \textcolor{Tan}{\textbf{beige}}, the Moon in \textcolor{Grey}{\textbf{grey}}, and the NRHO in \textcolor{NavyBlue}{\textbf{dark blue}}.
During the MTE, the spacecraft is forced to coast ballistically, shown with a dashed \textcolor{DodgerBlue}{\textbf{light blue}} line.
Figure \ref{fig: motivating_mass_throttle} displays the corresponding wet mass and throttle profiles for the post-spiral phase in both solutions (nominal profiles in \textcolor{Chartreuse}{\textbf{green}}; reoptimized profiles in \textcolor{Fuchsia}{\textbf{pink}}), with dashed lines indicating average throttle levels in each case. 
The reoptimized trajectory deviates significantly from the nominal, initially trailing behind due to the MTE. 
It quickly closes the gap, ultimately completing the mission 5 days earlier than the nominal, but this, however, increases the propellant requirement by 8 kilograms.
The corresponding throttle profiles reveal that while both solutions follow similar patterns, the reoptimized solution exhibits a slightly higher average throttle, by approximately 10\%, in comparison to the nominal solution. 

Although a viable solution was found in this case, such recovery may not always be possible. 
The spacecraft's ability to recover from such engine failures depends on several factors, including its current state, the available on-board propellant, and the remaining mission time, underscoring the importance of designing reliable contingencies during the preliminary mission design phase. 
Doing so will allow us to not only build in the necessary margins but also ensure that an alternative solution is readily available to deploy in the event that such outages disrupt the execution of a scheduled burn.


In the most general case, a robust optimal control problem involves designing control strategies to ensure feasibility or near-optimality in the presence of uncertainties, which can be broadly classified as aleatoric or epistemic. 
While the robust optimal control problem is of interest to the broader controls community, it has gained particular traction in the aerospace domain, where aleatoric uncertainties (e.g., navigational errors), epistemic uncertainties (e.g., imperfect knowledge about system/spacecraft parameters), and maneuver execution errors can significantly impact mission success. 
In practice, robustness is typically evaluated a posteriori through navigation analysis during the preliminary mission design phase, followed by iterative adjustments to the nominal solution. 
Common methods to improve robustness include the addition of empirical margins, reductions in duty cycles, or incorporation of coasting arcs along sensitive solution segments. 
However, these adjustments, which are often decoupled from the optimization process, can lead to suboptimal solutions with overly conservative margins, reducing overall mission efficiency.

Existing approaches to robust optimal control generally address these uncertainties by either reformulating the problem to account for worst-case adversarial uncertainty realizations or by considering a range of possible outcomes, ensuring the solution remains robust across diverse scenarios. 
These methods, such as min-max optimization, and sensitivity-based control, have been widely employed for spacecraft trajectory design problems with each approach tailored to specific uncertainties classes. 
Alternative strategies include probabilistic methods, which directly incorporate parameterized uncertainty distributions into the optimization process, enabling a trade-off between performance and robustness, and have been successfully applied to a variety of problems.
For example, a stochastic differential dynamic programming framework has been applied to LT spacecraft trajectory design, using linear feedback control policies for correction maneuvers and unscented transform \cite{julier_new_2000, ross_unscented_2015} for computing expectation of random variables \cite{ozaki_stochastic_2018}. 
Tube stochastic optimal control extends this approach to problems with nonlinear constraints by sequentially approximating stochastic processes as Gaussian and introducing chance constraints for both open- and closed-loop control \cite{ozaki_tube_2020}. 
Chance-constrained stochastic optimal control has also been applied to problems with Gaussian uncertainties, leveraging convex optimization for impulsive transfers \cite{oguri_robust_2021} and primer vector theory for LT transfers \cite{oguri_stochastic_2022}.
Covariance control has been been explored in some recent studies, demonstrating its effectiveness in spacecraft guidance and control \cite{okamoto_optimal_2018, yi_nonlinear_2020}. 
Although these methods show promise, their applicability to more general problems is often constrained by their reliance on assumptions about the uncertainty distribution.
The belief optimal control paradigm offers a more general framework that accommodates various uncertainty classes without any assumptions about their underlying distributions \cite{greco_robust_2022}. 
However, as the authors note, scaling this approach to LT thrust problems presents significant challenges, including the higher dimensionality of the optimization problem, and the reliance on local solvers, which struggle in multimodal and high-dimensional spaces, as are typically characteristic of robust LT trajectory design problems.

Within this broader problem domain, the problem of designing trajectories robust against maneuver execution errors constitutes a specific and highly relevant application of robust optimal control theory. 
MTEs, which occur due to operational disruptions, lead to significant deviations from nominal LT trajectories and pose a critical challenge for mission planning and execution. 
Addressing this issue requires incorporating MTEs as a unique class of uncertainty in the optimization process.
In this study, we focus on this class of problems, where the term `robust problems' are used to specifically describe problems addressing maneuver execution errors in LT missions, particularly MTEs.
Despite growing interest from both industry and academia in recent years, theoretical and algorithmic developments in robust LT trajectory design (i.e., the process of designing LT transfers resilient to MTEs), have been quite limited.
Based on their respective approach, existing literature in this area can be predominantly classified into two broad categories, namely, \emph{missed thrust analysis} and \emph{missed thrust design}. 
The missed thrust analysis approach focuses on evaluating a nominal solution's sensitivity to MTEs, often to diagnose potential risks and/or incorporate margin allocations, and in contrast, the missed thrust design approach focuses on designing the nominal solution by explicitly incorporating robustness as a performance metric during the optimization process.

Practical approaches to robust LT trajectory design are more similar to the former.
To account for MTEs, such approaches typically involve redesigning nominal trajectories, through lower duty cycles or coast arcs at strategically chosen points, and evaluating changes in key performance metrics, such as propellant usage and time of flight, to build in empirical margins (e.g., Dawn \cite{rayman_coupling_2007, oh_analysis_2008}, Psyche \cite{madni_missed_2020}).
Following a similar approach, Laipert et al. investigated the effects of single \cite{laipert_automated_2015} and multiple MTEs \cite{laipert_monte_2018} using Monte Carlo simulations.
However, the inherent decoupling between trajectory optimization and uncertainty quantification in these approaches risks transferring sensitivity to a different location along the trajectory, potentially resulting in suboptimal solutions.

In contrast, the missed thrust design approach involves directly integrating a robustness metric into the optimization framework, utilizing either deterministic or stochastic methodologies.
Deterministic approaches involve incorporating constraints on missed thrust recovery margins (i.e., the maximum amount of time a spacecraft may be allowed to coast while still being able to reach the terminal manifold once thruster operations are resumed), into the optimization process.
State-of-the-art approaches extend this concept by lifting the original optimal control problem to a higher dimensional space to solve for a \textit{reference} trajectory (the path we plan to fly) simultaneously with multiple \textit{realization} trajectories (the path we may be forced to fly should an MTE occur) from a-priori chosen points along the reference trajectory \cite{mccarty_missed_2020, venigalla_low_2020}. 
While these methods seem effective, they face computational challenges as the number of realizations increases. 
Adaptive algorithms help mitigate these challenges, but application to chaotic multibody gravitational systems remains difficult.
Stochastic approaches, despite being relatively understudied in this context, have been applied to missed thrust design as well.
For example, Olympio approached the missed thrust design problem by formulating it as a two-level stochastic optimal control problem \cite{olympio_designing_2010}, whereas Rubinsztejn et al. used the expected thrust function under the assumption of a known MTE distribution \cite{rubinsztejn_designing_2021}.
More recent approaches have looked into using data-driven methods in an attempt to learn the mapping between the spacecraft state after an MTE has occurred and the optimal control sequence going forward using neural networks \cite{rubinsztejn_neural_2020, izzo_real-time_2021} and reinforcement learning \cite{miller_low-thrust_2019, zavoli_reinforcement_2021}. 
These methods, however, only solve a local problem (i.e., small perturbations from the nominal), and are often limited in their ability to generalize to more complex gravitational environments where the LT trajectories are more sensitive to perturbations.

An alternative approach for missed thrust design could be to leverage the relationship between the robust LT trajectories and the natural dynamical flow in multibody dynamical systems. 
Building on this concept, Alizadeh and Villac investigated incorporating a penalty term in the objective function to minimize deviations from the natural dynamical flow \cite{alizadeh_sensitivity_2013}. 
However, the penalty term accounts for the cumulative deviation over the entire trajectory, potentially allowing some trajectory segments to diverge significantly from the natural dynamics. 
Despite being a promising research avenue, the exact relationship between the dynamical flow and resulting trajectories remained unclear, highlighting a gap in our understanding of robust LT trajectory design from a dynamical systems perspective.
This gap has been addressed in a recent study by Sinha and Beeson, which quantitatively demonstrated how robust LT solutions exploit dynamical structures in multibody systems, providing a foundation for optimization algorithms to leverage this information, and enabling more effective global search strategies for robust LT trajectory design \cite{sinha_statistical_2024_arxiv}.

Even in the absence of robustness considerations, the evolving operational requirements during the preliminary mission design phase necessitate rapid and efficient exploration of the solution space.
This exploration typically requires a global search algorithm, where generating high-quality solutions critically depends on having effective initial guesses. 
In an ideal scenario, a global search algorithm would place initial guesses within the basins of attraction for a collection of local optima, so that gradient-based solvers can converge reliably to qualitatively diverse solutions. 
However, identifying an optimal sampling distribution for these initial guesses is highly non-intuitive, particularly in high-dimensional LT trajectory design problems, and often problem-dependent, especially when the objective landscape evolves alongside the mission parameters in early design stages.

Global search problems in astrodynamics are typically solved using well-established algorithms, such as genetic algorithms, particle swarm optimization, differential evolution, simulated annealing, and basin-hopping methods \cite{locatelli_global_2013}. 
Significant efforts have been devoted to Pareto front discovery \cite{hartmann_optimal_1998, russell_primer_2007, oshima_global_2017} and search space analysis \cite{dilizia_advanced_2004, vasile_preliminary_2006, vasile_analysis_2010, addis_global_2011}, providing valuable insights into current global search challenges. 
For a comprehensive discussion on global search algorithms, in particular their application in astrodynamics problems, we direct the reader to the excellent topical discussion by Beeson et al. \cite{beeson_global_2024}.
Monotonic basin hopping, introduced by Wales and Doye \cite{wales_global_1997} and later specialized by Leary \cite{leary_global_2000}, is a widely adopted global search algorithm for problems in astrodynamics. 
Several efforts have been made to optimize this method for trajectory design problems. 
For instance, Englander and Englander have attempted to identify a heuristically optimal distribution for a given problem, where it was shown that a Pareto distribution outperforms Cauchy, Gaussian, and uniform sampling, though it required careful manual parameter tuning \cite{englander_tuning_2014}. 
Subsequent advancements addressed this limitation by introducing adaptive local hop distributions, which dynamically adjust the sampling parameters during the search process \cite{englander_automated_2017, englander_hopping_2020}.
It is well-established in existing literature that the effectiveness of global searches relies heavily on understanding the topology of the solution space, a challenge that becomes particularly pronounced in high-dimensional problems.
Additionally, as the problem dimensionality increases, the optimal basins shrink relative to the parameter space, making na\"ive sampling strategies insufficient. 
This challenge is further exacerbated in robust LT trajectory design, where the solutions inherently resides in a higher-dimensional space, increasing the complexity of identifying the local optimal basins.

To summarize, existing literature in missed thrust design demonstrates significant advancements in this research area, spanning practical margin-based methods, deterministic and stochastic frameworks, data-driven approaches, and dynamical systems-inspired techniques.
Despite these advancements, significant challenges remain in understanding and improving global search algorithms for missed thrust design, particularly in devising strategies to generate effective initial guesses.


In this study, we present a novel initial guess generation strategy to improve the global search for robust LT solutions in complex multibody dynamical environments. 
We summarize the main contributions below:
\begin{enumerate}[label=\arabic*.]
    \item We compare two initial guess generation strategies for global search problems in robust LT trajectory design: 
    \begin{enumerate}[label*=\arabic*.]
        \item a \emph{non-conditional global} strategy which uses random sampling from a fixed distribution with global support, enabling a broad exploration of the design space, and
        \item a \emph{conditional global} strategy that leverages previously solved non-robust or simpler partially robust solutions (i.e., with lower depths of robustness) to refine and narrow the exploration of the design space, with the aim of accelerating the generation of high quality robust solutions.
        Additionally, we present a comprehensive discussion on various mapping strategies for generating initial guesses through the conditional approach.
    \end{enumerate}
    \item We validate the initial guess generation strategies on a realistic LT cislunar transfer within a medium-fidelity dynamical model, focusing on robust, minimum-fuel transfers for the Power and Propulsion Element. 
    We highlight the inherent complexities in designing LT trajectories in the cislunar space, which becomes even more challenging in the presence of MTEs.
    \item We conduct a comprehensive statistical analysis to evaluate the performance of both methodologies across multiple algorithmic metrics, such as feasibility ratios, solving time, and solution quality (through objective function value).
    Our results demonstrate that initializing the global search procedure with the conditional global strategy using previously solved non-robust solutions yields statistically significant improvements in both computational efficiency and solution quality. 
\end{enumerate}


The paper is organized as follows.
In \S \ref{section: robust problem formulation}, we present the robust LT trajectory design framework we utilize in this study, and in \S \ref{section: initial guess generation strategies}, we formalize the two global search approaches central to this paper.
Then, we present the dynamical model in \S \ref{section: dynamical model}, and introduce the case study in \S \ref{section: problem setup}.
We statistically compare the global search approaches using key algorithmic performance metrics in \S \ref{section: results and discussion}.
Finally, we highlight the importance of this work, and discuss the limitations of the current approach in \S \ref{section: conclusion}.
\section{Robust Problem Formulation}
\label{section: robust problem formulation}

A succinct mathematical derivation of the MTE problem is provided in this section. 
The derivation here is a restriction of the more general robust optimal control problem where randomness and stochasticity can enter at many levels. 
The reader is encouraged to review Sinha and Beeson \cite{sinha_statistical_2024_arxiv}, which provides proper mathematical context for the missed thrust design problem within the general robust optimal control problem. 
For the reader's convenience, we repeat some of the mathematical formulation from Sinha and Beeson \cite{sinha_statistical_2024_arxiv} and make use of the same mathematical notation for consistency.  
In Sec. \ref{subsection: robust problem formulation: the robust (MTE) formulation} we provide the mathematical definition of the general robust MTE problem where a MTE may occur at any time along the reference trajectory, as well as the possibility for multiple MTEs. 
In Sec. \ref{subsection: robust problem formulation: the restricted robust (MTE) formulation} we restrict to the case where at most a finite number of times along the reference trajectory may incur an MTE and further restrict to the scenario where at most one MTE may occur for the mission. 
The justification of this restriction based on real mission data is also provided in this section. 
Section \ref{subsection: problem formulation: transcription into a nonlinear program} then provides complete details on the control transcription and numerical solver setup for the restricted finite realization case that will be studied in the results section, \S \ref{section: results and discussion}.

\subsection{The General MTE Formulation}
\label{subsection: robust problem formulation: the robust (MTE) formulation}

Let $(\Omega, \mF, \P)$ be a probability space and $\omega \in \Omega$ a random sample.
We aim to find an extremal control solution $\bm{u}^* \in \mathcal{U}^\Omega$, with $\mathcal{U}^\Omega$ an admissible control set, to minimize the Bolza-type cost functional,
\begin{align}
\label{eq: robust MTD}
\min_{\bm{u} \in \mathcal{U}^\Omega} \{ J(\bm{u}^\dagger) \equiv \phi(\bm{\xi}^\dagger_1) + & \int_0^1 \mathcal{L}(s, \bm{\bm{\xi}^\dagger_s}, \bm{u}^\dagger_s) ds \},
\end{align}
such that the dynamical constraints with deterministic boundary conditions, $\Xi_0, \Xi_1$, given by Eq. \eqref{eq: robust MTD dynamics equation} is satisfied, 
\begin{align}
\label{eq: robust MTD dynamics equation}
\bm{\xi}^\omega_t = \bm{\xi}^\dagger_0 + \int_0^t \bm{f}(s, \bm{\xi}^\omega_s) ds &+ \int_0^{\tau_1(\omega) \wedge t} \bm{g}(s, \bm{\xi}^\omega_s, \bm{u}^\dagger_s) ds \nonumber \\
&+ \sum_{i \in \integers_+} \int_{\tau_{2i}(\omega)}^{\tau_{2i + 1}(\omega) \wedge t} \bm{g}(s, \bm{\xi}^\omega_s, \bm{u}^\omega_s) ds, \quad \forall t \in [0, 1], \quad \forall \omega \in \Omega, \\
& \xi^\omega_0 \in \Xi_0, \quad \xi^\omega_1 \in \Xi_1. \nonumber
\end{align}
In Eq. \ref{eq: robust MTD}, $\bm{u}^\dagger$ is the \emph{reference} (or nominal) control solution and $\bm{\xi}^\dagger$ the reference state solution generated by $\bm{u}^\dagger$ that a mission designer would hope to fly. 
Although the reference solution determines the objective value of the problem, where $\mathcal{L}$ denotes the running cost and $\phi$ the terminal cost, it is coupled to a set of \emph{realization} solutions $\xi^\omega$ that must satisfy the same boundary conditions as $\bm{\xi}^\dagger$, as shown in Eq. \ref{eq: robust MTD dynamics equation}.
The realization solutions are the trajectories that would result from an MTE. 
The $f$ coefficient in Eq. \ref{eq: robust MTD dynamics equation} provides the vector field for the natural dynamics of the system, whereas the $g$ coefficient represents a vector field dependent on the control input. 
Lastly, the optimal control problem is normalized in time on the unit interval $[0, 1]$ for simplicity of presentation. 

All randomness for this optimal control problem enters via the random variable $\tau : \Omega \rightarrow [0, 1]$ in Eq. \ref{eq: robust MTD dynamics equation}.
In particular, this variable is defined as strictly increasing random times $\tau \equiv \{\tau_i(\omega) \in \mathbb{R}_+ \ | \ \tau_i < \tau_{i+1}, \ \forall i \in \mathbb{Z}_+, \ \omega \in \Omega \}$, where $\mathbb{Z}_+$ are the positive integers.
In what follows, we will also identify $\Omega$ with the unit circle (i.e., $\Omega \simeq S^1 \simeq [0, 1]$).
The symbol $\wedge$ in Eq. \eqref{eq: robust MTD dynamics equation} is the minimum operator (i.e., $a \wedge b = \min(a, b)$). 
Hence, for a given sample $\omega \in \Omega$, an MTE will be initiated if $\tau_1(\omega) < 1$, with additional MTEs for any $\tau_{2i + 1}(\omega) < 1$ with $i \in \mathbb{Z}_+$. 
The duration of an MTE is $\tau_{2(i + 1)} - \tau_{2i + 1}$ for any $i \in \mathbb{Z}_+ \cup \{0\}$. 

Our choice of admissible control sets start by defining the set $\mathcal{U} \equiv PC([0, 1]; \reals^n)$ for $n \in \integers_+$ to be the piecewise continuous functions on $[0, 1]$. 
Our admissible control set will then be given by $\mathcal{U}^\Omega \equiv \mathcal{U}^{S^1}$, which describes the functions from $S^1$ into $\mathcal{U}$ (equivalently $\prod_{S^1} \mathcal{U}$).   
We make the choice that $\tau_1(0) = \tau_1(1) > 1$, and hence the sample $\omega \in \{ 0, 1 \}$ will correspond to a (deterministic) non-MTE trajectory for Eq. \eqref{eq: robust MTD dynamics equation}.
We denote this special case, when $\omega \in \{0, 1\}$, with the $\dagger$ symbol as $\bm{\bm{u}^\dagger}$ and refer to the state solution $\bm{\bm{\xi}^\dagger}$ as the \emph{reference} solution. 
For all other cases, when $\omega \in (0, 1)$, we denote the control solution as $\bm{u^\omega}$ and refer to the associated state solution $\bm{\xi^\omega}$ as a \emph{realization} solution. 

\subsection{The Restricted (Finite Realization) MTE Formulation}
\label{subsection: robust problem formulation: the restricted robust (MTE) formulation}

As discussed in Sec. \ref{subsection: robust problem formulation: the robust (MTE) formulation}, the random time $\tau$, may characterize an arbitrary number of MTEs. 
The description of how many MTEs may occur is encoded in the definition of the probability distribution $\P$, or likewise the one that is induced on $[0, 1]$ by the random variable $\tau$. 
Based on analysis of past LT missions, Imken et al. \cite{imken_modeling_2018} have suggested that the Weibull distribution is a good fit for the initiation and duration times of a MTE. 
Because the Weibull distribution is a continuous distribution, achieving numerical tractability would require some sample approximation. 
In this paper, we adopt a simpler distribution, and our assumptions align with those by McCarty and Grebow \cite{mccarty_missed_2020}, as well as Venigalla et al. \cite{venigalla_low_2020}. 
In particular, we study the case where at most one MTE occurs during the mission, and the time of its initiation can be sufficiently approximated at three distinct points.
Further justification for these restrictions is given at the end of this section. 
First we provide all of our assumptions conceptually in words, followed by their corresponding mathematical formulations as subitems:

\begin{enumerate}[label=A\arabic*.]
\item One MTE will occur for any realization.
\label{assumption: A1}
\begin{itemize}
\item 
For each $\omega \in \Omega$, assume that $\tau_3(\omega) > 1$.
\end{itemize}
\item 
At most, three MTE initiations are considered, with each corresponding to the start of a thrust segment (a shooting transcription is used and will be further explained in \S \ref{subsection: problem formulation: transcription into a nonlinear program}).
\label{assumption: A2}
\begin{itemize}
\item 
Assume that $(0, 1) \subset S^1 = \Omega$ is partitioned into a collection of $N$ intervals $(E_i)_{i=1}^N$. 
\item 
Assume that for every interval $E_i$, that for any $\omega_0, \omega_1, \omega_2, \omega_3 \in E_i$, if $\tau_1(\omega_0) \neq \tau_1(\omega_1) \neq \tau_1(\omega_2)$, then we must have $\tau_1(\omega_3) = \tau_1(\omega_j)$ for at least one of $j \in \{0, 1, 2\}$.
\end{itemize}
\item Only a finite number of MTE durations are allowed.
\label{assumption: A3}
\begin{itemize}
\item Assume that each interval $E_i$ is further partitioned into a collection of $M$ subintervals $(E_{i, j})_{j = 1}^M$. 
\item Assume that for every subinterval $E_{i, j}$, that we have for any $\omega_0, \omega_1 \in E_{i, j}$, the relation $\tau_2(\omega_0) - \tau_1(\omega_0) = \tau_2(\omega_1) - \tau_1(\omega_1)$.
\end{itemize}
\item Lastly, for the $NM$ versions of the robust MTE optimal control problem just defined, we consider an equal number of different probability distributions $(\mathbb{P}_{i, j})$ on $\Omega$, with each having support only on their corresponding $E_{i, j}$ interval.  
\label{assumption: A4}
\begin{itemize}
\item For each $i \in \{1, \hdots, N\}, j \in \{1, \hdots, M\}$, define a robust MTE optimal control problem where $\mathbb{P}(E_{i, j} \cup \{0, 1\}) = 1$.
\end{itemize}
\end{enumerate}


We now provide justifications for the governing assumptions underlying our analysis.
Our first assumption \ref{assumption: A1} restricts our analysis to a single MTE per realization. 
While this may initially seem limiting, it is well supported by statistical evidence from historical missions \cite{imken_modeling_2018}. 
For instance, an analysis of the baseline solution for our case study (discussed in more detail in \S \ref{section: problem setup}) reveals that a single MTE accounts for approximately 90\% of all possible missed thrust scenarios.
The baseline solution comprises a maximum contiguous thrust arc lasting approximately 276 days, spanning both the spiral and post-spiral phases. 
Of these, the spiral phase lasts 178 days, leaving 98 days for the post-spiral phase. 
We focus solely on the post-spiral phase for this study, as this phase is more susceptible to perturbative effects from multibody dynamics, and assume that no MTEs occur during the spiral phase. 
According to the Weibull distribution parameters from Imken et al. \cite{imken_modeling_2018}, a single MTE encompasses approximately 90\% of possible MTE scenarios within a contiguous thrust arc lasting 98 days, justifying our assumption.

\ifthenelse{\boolean{includefigures}}
{
    \begin{figure}[!htb]
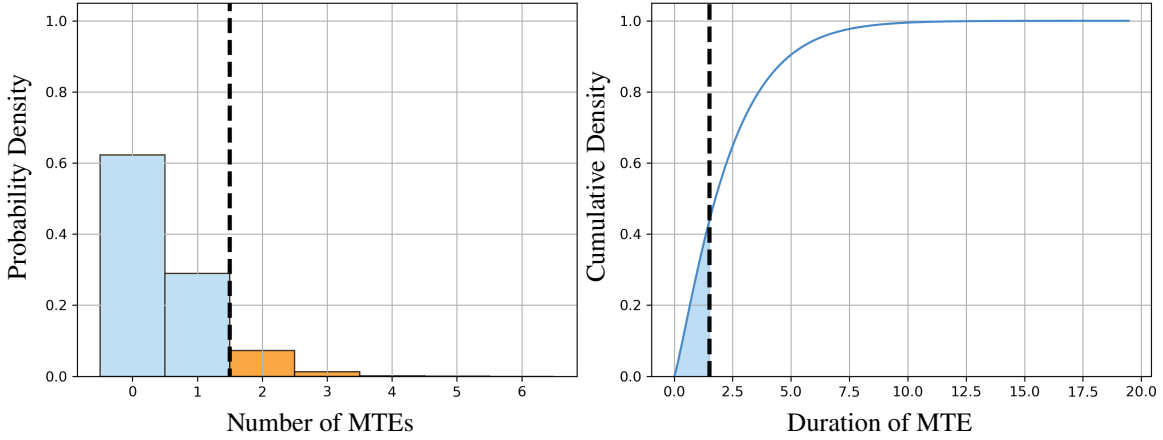

        \centering
        \begin{tikzonimage}[keepaspectratio, width=0.9\textwidth]{"distribution.png"}
        \node[fill=white, opacity=1.0, text opacity=1, anchor=south] at (0.25,-0.1) {Number of MTEs};
        \node[fill=white, opacity=1.0, text opacity=1, anchor=south, rotate=90] at (0.0,0.55) {Probability Density};
        \node[fill=white, opacity=1.0, text opacity=1, anchor=south] at (0.75,-0.1) {Duration of MTE};
        \node[fill=white, opacity=1.0, text opacity=1, anchor=south, rotate=90] at (0.5175,0.55) {Cumulative Density};
        \end{tikzonimage}
        \caption{Distribution of number of MTEs and their durations based on historical data (see Imken et al. \cite{imken_modeling_2018}).}
        \label{fig: distribution}
    \end{figure}
}
{
}

While a single MTE could theoretically occur at any point along the trajectory, simulating every possible outage location with non-zero probability would require a potentially infinite number of realizations, creating an insurmountable computational burden. 
Further details on how the computational complexity scales with with number of MTE initiations is given in \S \ref{subsection: problem formulation: transcription into a nonlinear program}. 
To address this, we restrict the analysis to a maximum of three initiation points along the reference trajectory, as outlined in Assumption \ref{assumption: A2}. 
These initiation points are strategically chosen to examine sensitive regions of the trajectory, specifically the beginning, middle, and end of the transfer. 
As outlined in Assumption \ref{assumption: A3}, we simulate three distinct MTE durations for each initiation point, ranging up to 1.5 days, representing approximately 40\% of the possible outage durations observed in historical data.
Finally, under Assumption \ref{assumption: A4}, we examine specific combinations of the selected MTE initiation points and durations. 
Detailed descriptions of the cases analyzed in this study are provided in \S \ref{section: results and discussion}.
The range of MTE scenarios analyzed in this study is depicted in Fig. \ref{fig: distribution}, with the blue shaded regions, marked by dashed lines, indicating the intervals considered. 
These regions highlight the number of MTEs and their corresponding durations, providing a balance between computational tractability and adequate coverage of critical scenarios.
Future work may expand this analysis to include a broader range of initiation points and outage durations.


\ifthenelse{\boolean{includefigures}}
{
    \begin{figure}[!htb]
        \centering
        \begin{subfigure}[b]{0.45\textwidth}
            \centering
            \includegraphics[keepaspectratio, width=\textwidth]{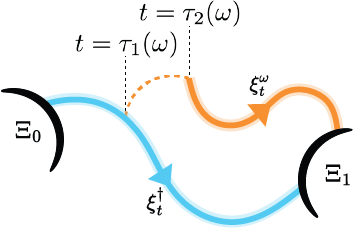}
            \caption{
            Schematic of the restricted robust problem showing \textcolor{SkyBlue}{reference} and \textcolor{Orange}{realization} trajectories
            }
            \label{fig: problem_setup_trajectory}
        \end{subfigure}
        \hfill
        \begin{subfigure}[b]{0.45\textwidth}
            \centering
            \includegraphics[keepaspectratio, width=\textwidth]{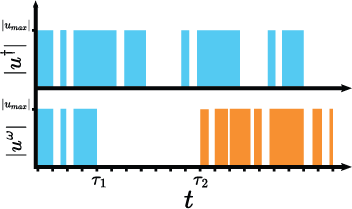}
            \caption{
            Schematic of corresponding throttle profiles with thruster outage during interval $t \in [\tau_1, \tau_2]$
            }
            \label{fig: problem_setup_throttle}
        \end{subfigure}
        \caption{
        Schematic of the restricted robust problem with one allowable MTE initiation location
        }
        \label{fig: problem_setup}
    \end{figure}
}
{
}

For completeness, we now state the optimal control problem for the restricted missed thrust design problem (see Fig. \ref{fig: problem_setup} for a schematic) under the additional assumptions just given:
\begin{align}
    \label{eq: simplified MTD}
    \min_{\bm{u} \in \mathcal{U}^\Omega} \{ J(\bm{u}^\dagger) \equiv \phi(\bm{\xi}^\dagger_1) + & \int_0^1 \mathcal{L}(s, \bm{\bm{\xi}^\dagger_s}, \bm{u}^\dagger_s) ds \},
\end{align}
where the reference and realization dynamics are given by,
\begin{align}
    \label{eq: simplified MTD dynamics equation}
    \quad \bm{\xi}^\omega_t = \bm{\xi}^\dagger_0 + \int_0^t \bm{f}(s, \bm{\xi}^\omega_s) ds &+ \int_0^{\tau_1(\omega) \wedge t} \bm{g}(s, \bm{\xi}^\omega_s, \bm{u}^\dagger_s) ds \nonumber \\
    &+ \int_{\tau_2(\omega)}^{t} \bm{g}(s, \bm{\xi}^\omega_s, \bm{u}^\omega_s) ds, \quad \forall \ t \in [0, 1], \quad \forall \omega \in \Omega, \\
    & \xi^\omega_0 \in \Xi_0, \quad \xi^\omega_1 \in \Xi_1. \nonumber
\end{align}

\subsection{Transcription into a nonlinear program}
\label{subsection: problem formulation: transcription into a nonlinear program}

To solve the missed thrust design problem, we make use of Dynamically Leveraged Automated (N) Multibody Trajectory Optimization (DyLAN), a computational astrodynamics software package developed by Beeson et al. \cite{beeson_dylan_2022}.
DyLAN brings together dynamical systems tools with local and global optimization methods to search for solutions of optimal control problems in multibody environments. 
We assume a finite-burn low-thrust model for the LT trajectory, wherein, the trajectory is divided into discrete segments, and a continuous thrust is applied over the duration of each segment during which both the magnitude and the direction remain constant.
With this model, we use a direct optimal control approach with a forward-backward multiple shooting algorithm to transform the optimal control problem in Eq. \eqref{eq: simplified MTD} into a nonlinear program (NLP).
The gradient-based numerical optimizer SNOPT \cite{gill_snopt_2005} is then used to solve the NLP using initial guesses generated through two different strategies (as detailed in \S \ref{section: initial guess generation strategies}).

We define the NLP as follows:
\begin{equation}
    \label{eq: nonlinear program}
    \begin{split}
    &\quad \min \limits_{\bm{x}^\dagger \in \mathbb{R}^{\overline{N}^\dagger}, \ \bm{x}^\omega \in \mathbb{R}^{\overline{N}^\omega}} \{ J(\bm{x}^\dagger) = -m_f^\dagger \}, \\
    \text{subject to} &\quad c^\dagger_l(\bm{x}^\dagger) = 0, \quad c^\omega_l(\bm{x}^\omega) = 0, \quad \forall \ l \in \mathcal{E}, \\
    &\quad c^\dagger_l(\bm{x}^\dagger) \leq 0, \quad c^\omega_l(\bm{x}^\omega) \leq 0, \quad \forall \ l \in \mathcal{I},
    \end{split}
\end{equation}
where the index set $\mathcal{E}$ identifies the equality constraints, consisting of matchpoint defect errors for the position, velocity, and mass continuity of the reference and realization solutions. 
The index set $\mathcal{I}$ identifies the inequality constraints, which consists of simple bounds on the components of the reference control decision variable $\bm{x}^\dagger$ and realization control decision variable $\bm{x}^\omega$. 
The reference control decision variable has $\overline{N}^\dagger = 3 N^\dagger + 4$ components given by, 
\begin{equation}
\label{equation: decision vector}
\bm{x}^\dagger = (T^\dagger_s, T^\dagger_i, T^\dagger_f, \bm{u}^\dagger_1, \bm{u}^\dagger_2,..., \bm{u}^\dagger_{N^\dagger}, m^\dagger_f),
\end{equation}
where $T^\dagger_s$ is the shooting time, $T^\dagger_i$ the initial coast time, $T^\dagger_f$ the final coast time, and therefore the total time-of-flight is $T^\dagger_i + T^\dagger_s + T^\dagger_f$. 
$N^\dagger$ represents the number of finite burn thrust segments.
$\bm{u}^\dagger_p \in \mathbb{R}^3$ is a constant thrust vector (characterized in spherical coordinates by the throttle, in-plane, and out-of-plane thrust angle) for the $p^\text{th}$ thrust segment where $p \in \{1, 2, \cdots, N^\dagger\}$. 
The thrust segments each have equal time of $T^\dagger_s / N^\dagger$. 
Lastly, $m^\dagger_f$ is the final delivered wet mass. 
The total number of constraints for the reference solution is equal to $\overline{N}^\dagger + 7$, comprising $N^\dagger$ simple bound constraints and 7 equality constraints.
A similar transcription follows for the non-robust control solution.

The transcription of the realization solutions also follows a similar approach but includes fewer control variables. 
we use an \emph{adaptive segmentation} strategy to discretize the realization control solutions, which dynamically adjusts the number of control segments depending on where the MTE begins, to ensure consistency in control authority between the reference and the realization solutions.
To illustrate, consider a random sample $\omega_k \in \Omega$ corresponding to a realization solution $\bm{\xi}^{\omega_k}$.
We assume that the realization solution $\bm{\xi}^{\omega_k}$ starts at the beginning of the $n_{\omega_k}^{\text{th}}$ control segment within the reference solution $\bm{\xi}^\dagger$, where $1 \leq n_{\omega_k} \leq N^\dagger$.
In other words, for this particular realization, $\tau_1(\omega_k)$ coincides with the beginning of the $n_{\omega_k}^{\text{th}}$ reference control segment, and the value of $\tau_1(\omega_k)$ determines the number of segments in $\bm{\xi}^{\omega_k}$ given by $\overline{N}^{\omega_k} = 3 N^{\omega_k} + 4 = 3(N^\dagger - n_{\omega_k}) + 4$.
Recall that $\bm{x}^{\omega_k}$ contains $N^{\omega_k}$ thrust segments, with each thrust segment having an equal duration of $T^{\omega_k}_s / N^{\omega_k}$. 
Since it is expected that $T^{\omega_k}_s \leq T^\dagger_s$, if we naively assign $N^{\omega_k} = N^\dagger$, it would imply that each realization control segment now may span a shorter time interval (since $T^{\omega_k}_s / N^\omega \leq T^\dagger_s / N^\dagger$) which would inadvertently result in a higher control authority for $\bm{x}^{\omega_k}$ compared to $\bm{x}^{\dagger}$.
So, instead, the number of segments in $\bm{x}^{\omega_k}$ is adjusted by setting $N^{\omega_k} = N^\dagger - n_{\omega_k}$, such that each thrust segment then spans comparable time intervals, and thereby help promote congruity in the control authority between the reference and the realization solutions.
The realization control decision variable has $\overline{N}^{\omega_k} = 3 N^{\omega_k} + 4 = 3 (N^\dagger - n_{\omega_k}) + 4$ components given by, 
\begin{equation}
\label{equation: realization decision vector}
\bm{x}^{\omega_k} = (T^{\omega_k}_s, T^{\omega_k}_i, T^{\omega_k}_f, \bm{u}^{\omega_k}_1, \bm{u}^{\omega_k}_2,..., \bm{u}^{\omega_k}_{N^{\omega_k}}, m^{\omega_k}_f),
\end{equation}

\begin{figure}[hbt!]
    \centering
    \includegraphics[width=0.75\linewidth]{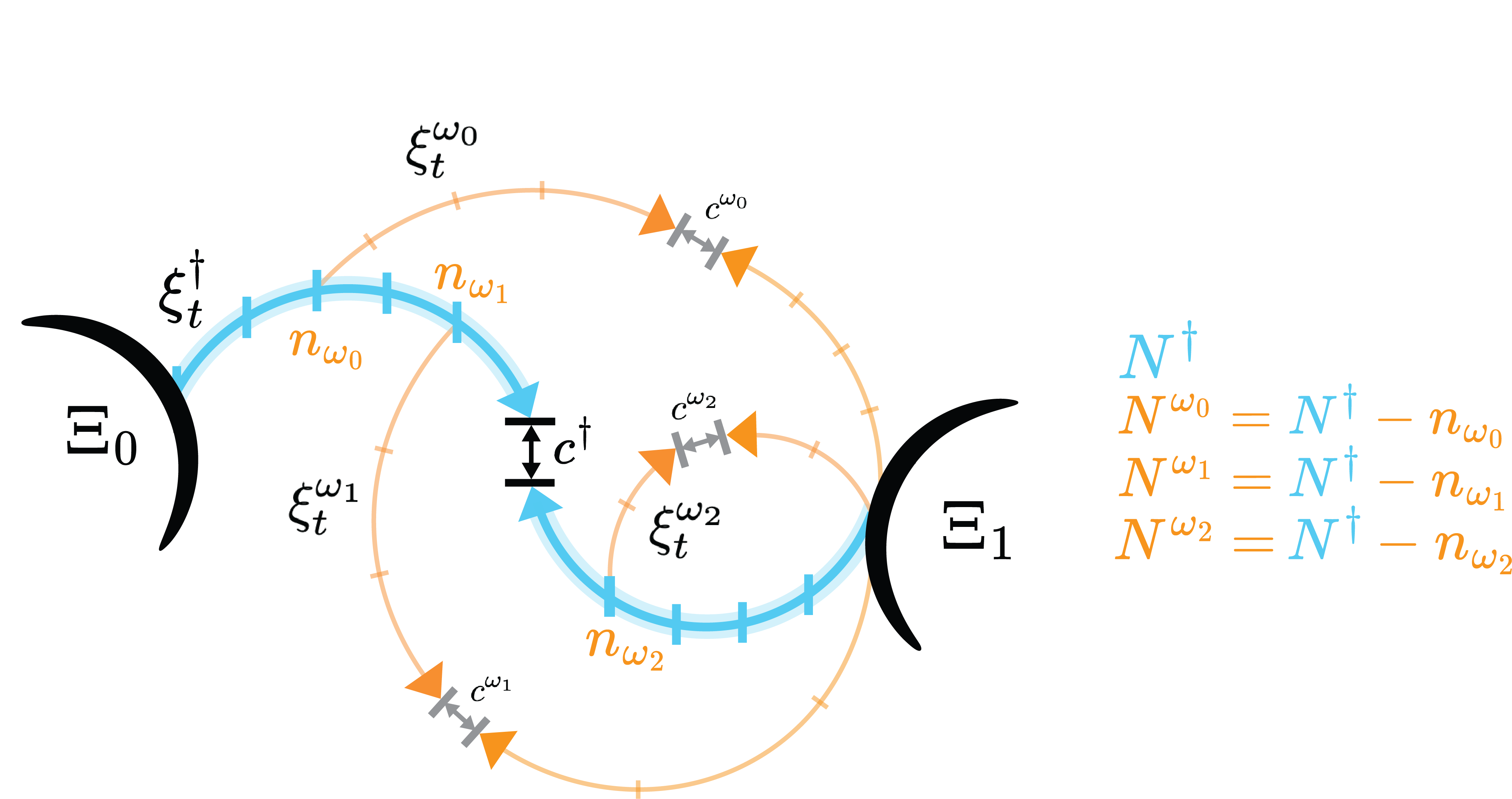}
    \caption{
    Dynamic allocation of realization segments using adaptive segmentation
    }
    \label{fig: adaptive_segmentation}
\end{figure}

The adaptive segmentation strategy is illustrated in Fig. \ref{fig: adaptive_segmentation} which demonstrates an example with three random samples, $\{\omega_0, \omega_1, \omega_2\} \in \Omega$, which yields the realization solutions $\{\bm{\xi}^{\omega_0}, \bm{\xi}^{\omega_1}, \bm{\xi}^{\omega_2}\}$. 
We assume a forward-backward shooting transcription for all reference and realization solutions, but the idea is applicable to forward shooting transcription as well.
In a forward-backward shooting phase, the spacecraft state is propagated from both boundaries inwards in contrast to a forward-shooting algorithm where the spacecraft state is propagated from one of the boundaries directly to the other.
In general, this two-sided shooting results in a discontinuity of the spacecraft’s state vector at some location along the phase, which we attempt to minimize using a numerical optimizer.
In Fig. \ref{fig: adaptive_segmentation}, the matchpoint defect for the reference solution is denoted by $\bm{c}^\dagger$, and by $\bm{c}^{\omega_k}$ for the $k$\textsuperscript{th} realization solution $\xi^{\omega_k}$. 
The reference solution $\bm{x}^\dagger$ is discretized into $N^\dagger$ segments. For the realization solution $\bm{\xi}^{\omega_0}$, which begins at the $n_{\omega_0}^{\text{th}}$ segment, the number of segments $N^{\omega_0}$ is reduced by $n_{\omega_0}$. 
Similar adjustments are made for $\bm{\xi}^{\omega_1}$ and $\bm{\xi}^{\omega_2}$.
The full control decision variable for the problem has $\overline{N}^\dagger + \sum_{k=0}^{2} \overline{N}^{\omega_k}$ components given by,
\begin{align}
    \bm{x} = (T^\dagger_s, T^\dagger_i, T^\dagger_f, \bm{u}^\dagger_1, \bm{u}^\dagger_2,..., \bm{u}^\dagger_{N^\dagger}, m^\dagger_f, & T^{\omega_0}_s, T^{\omega_0}_i, T^{\omega_0}_f, \bm{u}^{\omega_0}_1, u^{\omega_0}_2,..., \bm{u}^{\omega_0}_{N^{\omega_0}}, m^{\omega_0}_f, \nonumber \\ 
    & T^{\omega_1}_s, T^{\omega_1}_i, T^{\omega_1}_f, \bm{u}^{\omega_1}_1, u^{\omega_1}_2,..., \bm{u}^{\omega_1}_{N^{\omega_1}}, m^{\omega_1}_f, \nonumber \\ 
    & T^{\omega_2}_s, T^{\omega_2}_i, T^{\omega_2}_f, \bm{u}^{\omega_2}_1, u^{\omega_2}_2,..., \bm{u}^{\omega_2}_{N^{\omega_2}}, m^{\omega_2}_f),
\end{align}

A breakdown of the number of decision variables for both the non-robust and robust problems (with $K$ realizations for generality) are summarized in Table \ref{tab: decision_variables}. 

\begin{table}[hbt!]
    \caption{Number of Decision Variables (Number of Realizations = K)}
    \label{tab: decision_variables}
    \centering
    \begin{tabular}{lcc}
        \hline
        & Non-Robust & Robust \\\hline
        Number of Segments & $N^\dagger$ & $N^\dagger$+$\sum_{k=0}^{K - 1} N^{\omega_k}$ \\
        \multicolumn{3}{l}{Control Vector Components} \\
        \hspace{3mm} Time of Flight & 3 & 3 + 3K \\
        \hspace{3mm} Thrust Vector & 3$N^\dagger$ & 3($N^\dagger$+$\sum_{k=0}^{K - 1} N^{\omega_k}$) \\
        \hspace{3mm} Final Mass & 1 & 1 + K \\
        Number of Constraints & 7 & 7 + 7K \\
        \hline
    \end{tabular}
\end{table}

To solve the NLP in Eq. \ref{eq: nonlinear program}, it is necessary to compute the derivatives of the matchpoint defect errors with respect to the control variables. 
These derivatives can be determined either analytically or using approximation techniques, such as finite differences. 
Selecting an appropriate step size for the finite differences method, which balances numerical accuracy and solution time, is challenging and can lead to longer convergence times.
Analytic derivatives, if available, can significantly improve the efficiency of the optimization process and reduce computational overhead. 
In this study, analytic derivatives are provided to the optimizer for both non-robust and robust solutions. 
In the context of non-robust LT trajectory design, analytic derivatives have been extensively studied by Ellison et al. \cite{ellison_robust_2018}.
To compute the analytic derivatives for the non-robust problem, we adopt a similar methodology to their approach.
However, for robust LT trajectory design, it becomes necessary to augment the dynamical information by incorporating the flow of derivatives from the reference solution to the realization solutions.
For further details on computing the analytic derivatives for the robust problem, we refer the reader to Sinha and Beeson \cite{amlans_derivatives_2024}.

\ifthenelse{\boolean{includefigures}}
{
    \begin{figure}[!htb]
        \centering
        \begin{tikzonimage}[keepaspectratio, width=0.95\linewidth]{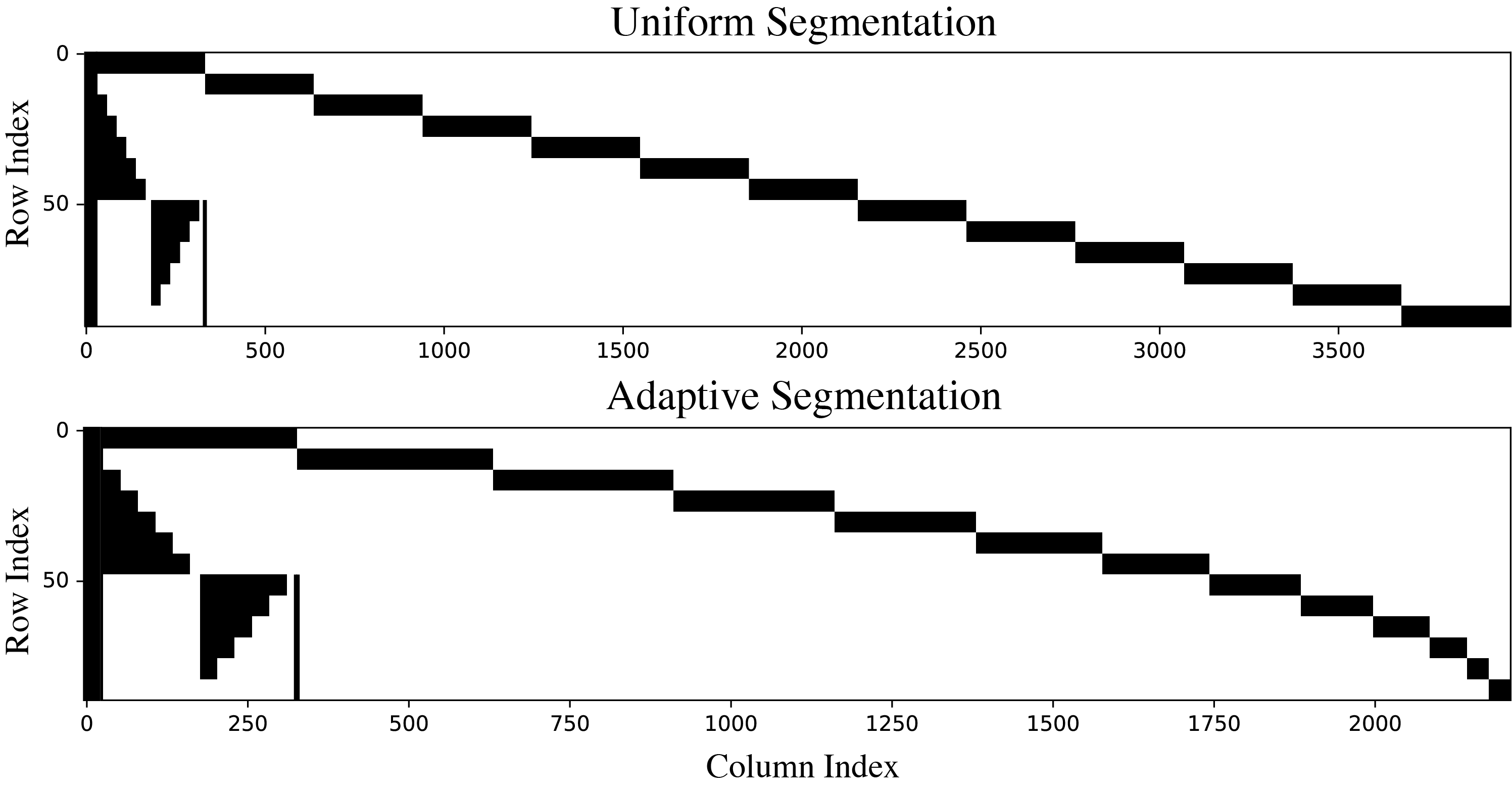} 
        \end{tikzonimage}
        \caption{Comparative analysis of Jacobian sparsity structure between adaptive and uniform segmentation}
        \label{fig: analytic derivatives sparsity pattern}
    \end{figure}
}
{
}

The matrix containing the partial derivatives of the matchpoint defect errors with respect to the control decision variables is known as the \emph{Jacobian}.
An additional benefit of the adaptive segmentation strategy can be illustrated by visualizing the sparsity pattern of the Jacobian matrix. 
To illustrate, consider a representative low-thrust transfer with $N^\dagger = 100$. 
We assume there exists twelve realization solutions, denoted as $\bm{\xi}^\omega = \{\bm{\xi}^{\omega_0}, \bm{\xi}^{\omega_1}, \dots, \bm{\xi}^{\omega_{11}}\}$, where $\tau_k$ coincides with the beginning of the thrust segments bearing indices $n_{\omega_k} = 8 + 8 k$ for $k = 0, 1, \dots, 11$.
The resulting sparsity pattern of the Jacobian matrix reveals key computational characteristics, such as the number of nonzero entries required for matrix operations.
Figure \ref{fig: analytic derivatives sparsity pattern} shows a comparison of the Jacobian sparsity patterns obtained using the adaptive segmentation strategy and the uniform segmentation strategy (i.e., without adjustments to the realization control solution) for the same problem. 
Even in this relatively simple example, the number of nonzero entries is significantly large. 
While this may not pose a challenge for standard computational algorithms, it can substantially impact the runtime of an NLP solver, which depends on repetitive and efficient matrix operations.
As illustrated in the figure, the adaptive segmentation strategy helps mitigate this issue by significantly reducing the number of dense entries in the Jacobian matrix.

\ifthenelse{\boolean{includefigures}}
{
    \begin{figure}[!htb]
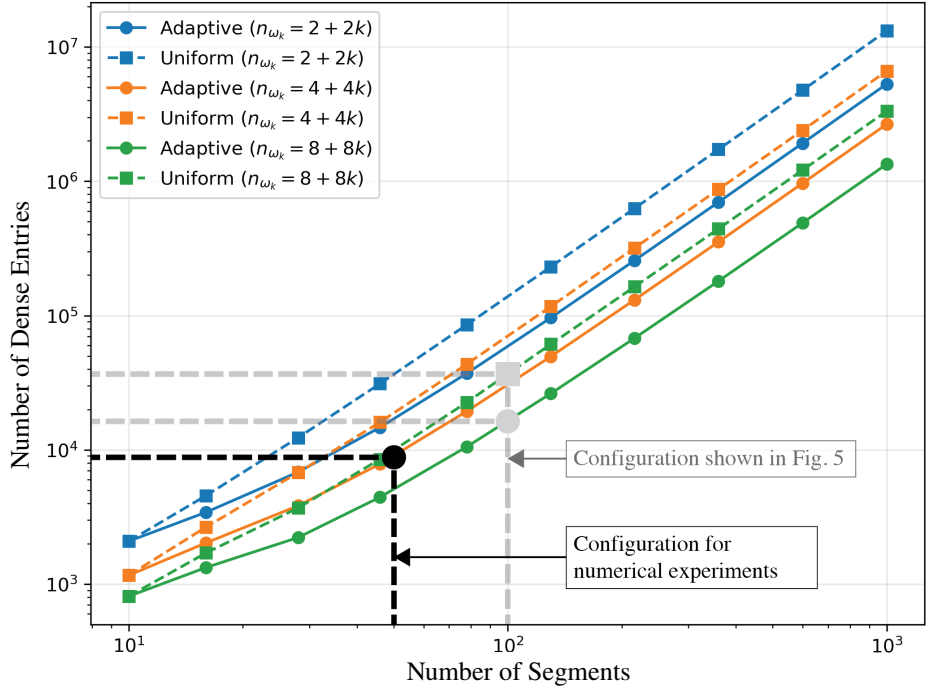

        \centering
        \begin{tikzonimage}[keepaspectratio, width=0.75\textwidth]{"number_of_dense.png"}
        \node[fill=white, opacity=1.0, text opacity=1, anchor=south] at (0.55,0.0) {Number of Segments};
        \node[fill=white, opacity=1.0, text opacity=1, anchor=south, rotate=90] at (0.035,0.51) {Number of Dense Entries};
        \end{tikzonimage}
        \caption{Comparative analysis of Jacobian density scaling between adaptive and uniform segmentation}
        \label{fig: dense entries}
    \end{figure}
}
{
}

To further understand the computational benefits of the adaptive segmentation strategy in addressing these challenges, let us now consider three different scenarios. 
Figure \ref{fig: dense entries} illustrates the relationship between the number of segments and the resulting number of dense entries in the Jacobian matrix, comparing adaptive and uniform segmentation approaches across three different scenarios.
The scenario labeled $n_{\omega_k} = M + Mk$ represents a configuration where the $k$\textsuperscript{th} realization initiates at the beginning of the $n_{\omega_k}$\textsuperscript{th} reference control segment, or in other words, realizations occur at every $M$ segments, resulting in $\left\lfloor \frac{N^\dagger}{M} \right\rfloor$ total realizations.
Results are shown for $M = \{2, 4, 8\}$.
The plot, presented on a logarithmic scale for both axes, demonstrates that the number of dense entries grows exponentially with increasing segment count. 
The adaptive segmentation strategy consistently yields fewer dense entries compared to its uniform counterpart across all configurations, with the disparity becoming more pronounced at higher segment counts. 
Two specific scenarios are highlighted: one used for subsequent numerical experiments in this study (indicated by the black dashed line) and another corresponding to the scenario in Fig. \ref{fig: analytic derivatives sparsity pattern} (shown by the gray dashed line). 
This analysis highlights that the adaptive approach offers substantial computational benefits by minimizing the density of the Jacobian matrix, a feature that becomes particularly advantageous as the complexity of the problem increases..

Another key observation from the sparsity pattern is the presence of off-diagonal terms (Fig. \ref{fig: analytic derivatives sparsity pattern}), which arise due to the coupling between realization defect constraints and reference decision variables. 
These terms reduce the sparsity of the Jacobian matrix, thereby diminishing the efficiency of sparse linear solvers. 
Numerical solvers perform optimally when the Jacobian exhibits a block-diagonal or nearly diagonal structure, as they can exploit sparsity for efficient linear matrix operations. 
The off-diagonal terms increase the density within the matrix blocks, leading to potentially higher space and time complexity for the solver, and thereby highlighting the additional numerical challenges inherent in addressing missed thrust design problems.
\section{Initial Guess Generation Strategies}
\label{section: initial guess generation strategies}


Efficient global search for robust solutions requires well-defined strategies for generating initial guesses that lead to a sufficient number of high-quality solutions. 
In this section, we explore two distinct methodologies, which we refer to as the \emph{non-conditional} and \emph{conditional} strategies, for constructing initial guesses tailored to the problem's requirements.
To define a mathematical notation for these strategies, we first define $\mathcal{P}_{k}$ as a robust problem with $0 \leq k \leq 3$ realizations. 
The case $k = 0$, refers to the non-robust (non-MTE) problem. 
Our two strategies are then defined as:
\begin{enumerate}
\item 
$\mathcal{S}(k)$, which represents the \emph{non-conditional} strategy applied to $\mathcal{P}_k$,
\item 
$\mathcal{S}(k \mid k')$, which represents the \emph{conditional} strategy applied to $\mathcal{P}_k$, where the initial guesses are conditioned on (informed by) solutions to $\mathcal{P}_{k'}$, with $0 \leq k' < k$.
\end{enumerate}

The defining difference between the non-conditional and conditional strategies is that the conditional strategy leverages data of already solved problems, and in particular, problems with fewer realizations and hence less depth of robustness being considered. 
Fundamentally, one would expect that there are tradeoffs of using the conditional approach in comparison to the non-conditional. 
For one, if we assume that the non-conditional approach is characterized by an initial guess generator $X : \Omega \rightarrow \mathcal{X}$, from the sample space $\Omega$ of the probability triple $(\Omega, \mathcal{F}, \mathbb{P})$ to the transcribed optimal control space $\mathcal{X}$ (some finite-dimensional Euclidean space), that induces a probability measure with full support $(\mathcal{X} = \operatorname{supp}\left( \mathbb{P}_X \equiv X_{*} \mathbb{P} \right))$, then this initial guess generator will eventually uncover all the optimal solutions if enough samples (initial guesses) are drawn. 
For problems with increasing numbers of realizations, the transcribed control space $\mathcal{X}$ will become larger and the optimal control problem more constrained due to the increasing coupling of the reference and realizations. 
Therefore, we expect that the non-conditional approach will require ever increasing number of samples to uncover a sufficient number of optimal solutions. 
The conditional strategy therefore makes use of known optimal solutions to simpler problem with less robustness (i.e., fewer realizations) to attempt to combat this curse of dimensionality and guide initial guess to the reduced space where feasible solutions exist. 
A potential drawback of the conditional approach, however, is that potential robust solutions may not be uncovered in the search process due to the narrower search domain. 

\subsection{Conditional and Non-Conditional Global Search}
\label{subsection: algorithmic considerations: conditional and non-conditional global search}

If $\mathcal{X}_k$ represents the transcription space for the $k$-th realization ($0$-th being the reference) and $\prod \mathcal{X}_k$ the space for the coupled reference and realizations, then the non-conditional strategy is characterized by a fixed distribution $\P_{X_k}$ on this product space, where $X_k : \Omega \rightarrow \prod \mathcal{X}_k$, with natural definition for ever increasing $k$. 
In the simplest form, $\P_{X_k}$ might just be the Uniform distribution, or other unimodal distributions such as the Pareto, Cauchy, or Gaussian which are often used with the popular Monotonic Basin Hopping (MBH) algorithm \cite{wales_global_1997, leary_global_2000, englander_automated_2012, englander_tuning_2014, englander_automated_2017}. 
A realization of our random variable $\P_{X_k} \sim \bm{x} = X_k(\omega)$ is an initial guess for an NLP solver, which we denote as $\pi : \mathcal{X}_k \rightarrow \mathcal{X}_k$, and is a map on the transcribed control space to produce optimal solutions (e.g., $\bm{x}^* = \pi(\bm{x})$). 
Composing $\pi$ with the random guess generation process, we get a random variable $X^*_k = \pi \circ X_k  : \Omega \rightarrow \mathcal{X}_k$, which itself induces a probability distribution $\P_{X^*_k}$.
This completes the non-conditional strategy.

Now, assume we have two problems $\mathcal{P}_{k'}$ and $\mathcal{P}_{k}$ with $0 \leq k' < k$, and a mapping $\mathcal{M}_{k'}^k : \mathcal{X}_{k'} \rightarrow \mathcal{X}_k$.
Then our conditional global search can be defined as the construction of the random variable $X_k \equiv \mathcal{M}_{k'}^k \circ \pi \circ X_{k'} : \Omega \rightarrow \mathcal{X}_k$ or similarly, the \emph{conditional} probability distribution $\P(\mathcal{X}_k \mid X_{k'})$ on $\mathcal{X}_k$. 
It is often the case that $\P(\mathcal{X}_k \mid X_{k'})$ is a Dirac distribution.

Returning to our informal definitions of $\mathcal{S}(k)$ and $\mathcal{S}(k \mid k')$ at the start of this section, these are now expressed as, 
\begin{align}
\mathcal{S}(k) \equiv \pi \circ X_k , \quad \text{and} \quad \mathcal{S}(k \mid k') \equiv \pi \circ \mathcal{M}_{k'}^k.
\end{align}



\begin{table}[!htbp]
\centering
\caption{Control variables bounds for the reference solution $\bm{x}^\dagger$ and realization solution $\bm{x}^\omega$}
\label{tab: bounds_nonconditional_global_search}
\begin{tabular}{l ccccc ccccc}
\hline
            & \multicolumn{5}{c}{\textbf{Reference}}                                                                                                              & \multicolumn{5}{c}{\textbf{Realization}}                                                                                                                                            \\ \hline
            & $T^\dagger_s$               & $T^\dagger_i$               & $T^\dagger_f$               & $\bm{u}^\dagger_p$               & $m^\dagger_f$               & $T^\omega_s$                                               & $T^\omega_i$               & $T^\omega_f$               & $\bm{u}^\omega_p$               & $m^\omega_f$               \\ \hline
Lower Bound & 0.0                         & 0.0                         & 0.0                         & $\bm{u}^\dagger_{p, \text{min}}$ & 0.0                         & 0.0                                                         & 0.0                         & 0.0                         & $\bm{u}^\dagger_{p, \text{min}}$ & 0.0                         \\
Upper Bound & $T^\dagger_{s, \text{max}}$ & $T^\dagger_{i, \text{max}}$ & $T^\dagger_{f, \text{max}}$ & $\bm{u}^\dagger_{p, \text{max}}$ & $m^\dagger_{f, \text{max}}$ & $(\frac{N^\dagger- n_\omega}{N^\dagger}) T^\dagger_{s, \text{max}}$ & $T^\dagger_{i, \text{max}}$ & $T^\dagger_{f, \text{max}}$ & $\bm{u}^\dagger_{p, \text{max}}$ & $m^\dagger_{f, \text{max}}$ \\ \hline
\end{tabular}
\end{table}

\noindent
The non-conditional strategy $\mathcal{S}(k)$ relies on sampling the components of the control decision variables $\bm{x}$ from bounded uniform distributions defined over the space $\prod \mathcal{X}_k$. 
The bounds for the control decision variables used in subsequent numerical experiments in this study, applicable to both the non-conditional strategy $\mathcal{S}(k)$ and the conditional strategy $\mathcal{S}(k \mid k')$, are provided in Table \ref{tab: bounds_nonconditional_global_search}. 
For the reference control solution $\bm{u}^\dagger$ and the realization control solution $\bm{u}^\omega$, the bounds are set to span the full range of the admissible control space, $\mathcal{U}$. 
To enforce a shorter maximum allowable shooting time of the realization solution compared to the reference solution, the upper bound of $T_s^\omega$ is truncated by the ratio $(N^\dagger - n_\omega)/N^\dagger$, where $n_\omega$ denotes the reference segment index where the realization begins.
The lower bound remains unchanged. 
The bounds for remaining realization components are identical to those of their corresponding reference values.



\subsection{Mapping Strategies For Conditional Global Search}
\label{subsection: algorithmic considerations: conditional mapping strategies}

\ifthenelse{\boolean{includefigures}}
{
    \begin{figure}[!htb]
        \centering
        \begin{tikzonimage}[keepaspectratio, width=0.95\linewidth]{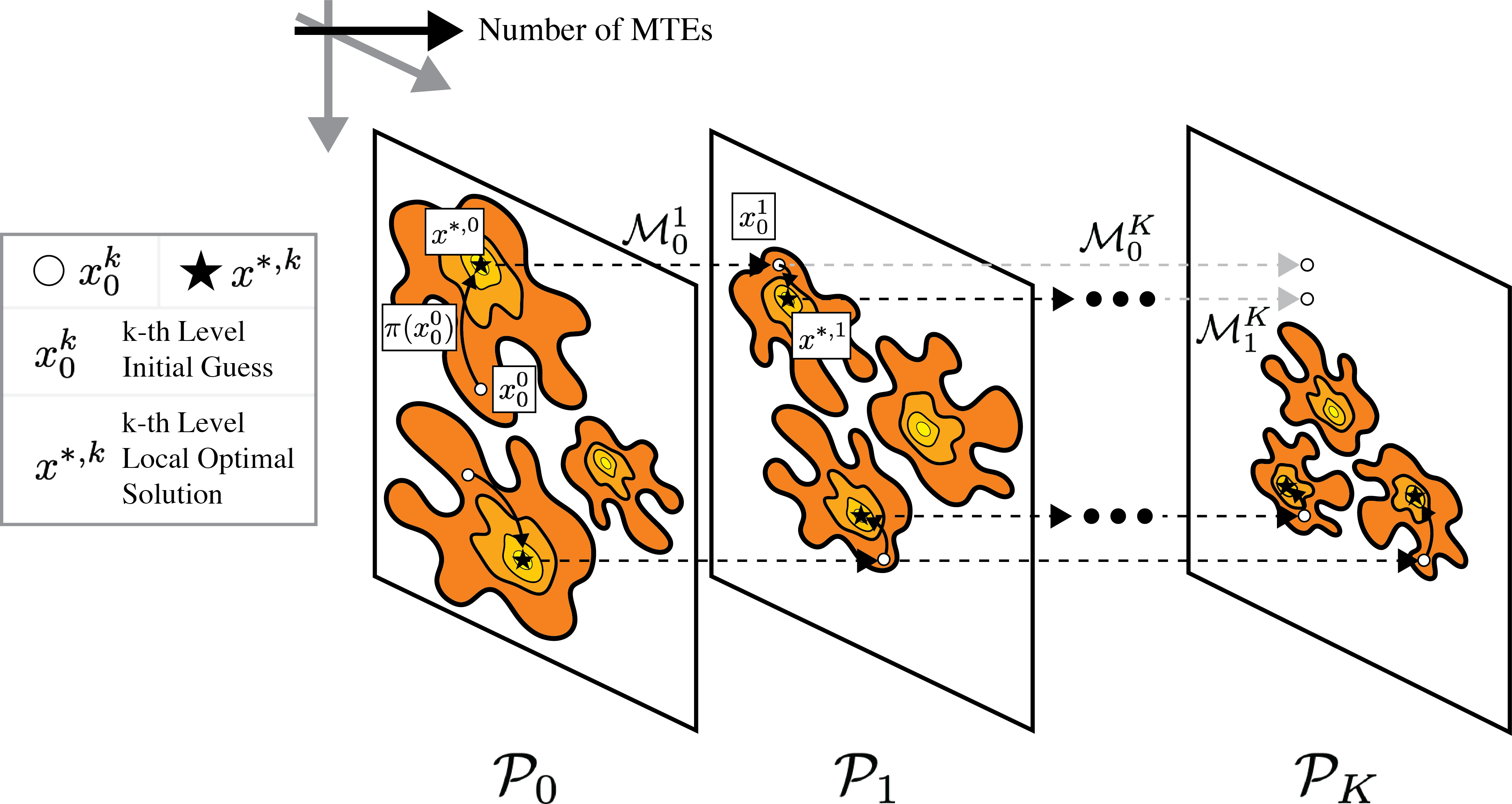} 
            \small
        \end{tikzonimage}
        \caption{Evolution of the objective function landscape with increasing depths of robustness}
        \label{fig: objective function landscape}
    \end{figure}
}
{
}

The conditional strategy is illustrated in Figure \ref{fig: objective function landscape}, which depicts the evolution of the objective function landscape as the problem transitions from a non-robust formulation ($\mathcal{P}_0$) to increasingly robust formulations ($\mathcal{P}_{k=1, \dots, K}$), each incorporating higher depths of robustness. 
Each vertical slice in the figure corresponds to the objective function landscape for a distinct problem $\mathcal{P}_k$.
The contours illustrate the objective function value, with the shaded areas depicting basins of attraction for different local minima, and the unshaded areas indicating infeasible regions in the solution space.
Initial guesses at each level $k$ are denoted by $\bigcirc$, while the optimal solutions to which they converge are represented by $\bigstar$. 
Arrows indicate the path taken by an initial guess $\bm{x}_0^k$, if initialized appropriately within a basin, as it evolves under the action of $\pi$. 

For problem $\mathcal{P}_0$, we sample an initial guess $\bm{x}_0^0 \in \mathcal{X}_0$ from the distribution $\mathbb{P}_{X_0}$. 
This point, residing within a basin of attraction, eventually converges to the local minimum $\bm{x}^{*, 0}$ under the solver $\pi$.
Problem $\mathcal{P}_1$ introduces a single realization solution, which mildly distorts the objective function landscape compared to $\mathcal{P}_0$. 
To initialize the solution process for $\mathcal{P}_1$, the optimal solution from the previous level, $\bm{x}^{*, 0}$, is projected onto the space $\mathcal{X}_1$ using a map $\mathcal{M}_0^1$, yielding the initial guess $\bm{x}_0^1 = \mathcal{M}_0^1(\bm{x}^{*, 0}) \in \mathcal{X}_1$.
Since this point also lies within a basin of attraction, it eventually converges to the local optimal solution $\bm{x}^{*, 1}$ under $\pi$. 
As we continue to increase $k$, and move on to more complex robust problems $\mathcal{P}_K$, the distortion of the objective function landscape can be so pronounced that it bears little resemblance to the original.
In this case, applying a similar map $\mathcal{M}_1^K$ and $\mathcal{M}_0^K$ to the optimal solutions $\bm{x}^{*, 1}$ and $\bm{x}^{*, 0}$ respectively, may result in initial guesses which now reside in infeasible regions of $\mathcal{X}_K$.
This phenomenon illustrates an inherent limitation with the conditional approach.
While it enables an efficient sequential algorithmic procedure for generating initial guesses, it cannot guarantee feasibility preservation as the depths of robustness increase.

\ifthenelse{\boolean{includefigures}}
{
    \begin{figure}[!htb]
        \centering
        \begin{tikzonimage}[keepaspectratio, width=0.75\linewidth]{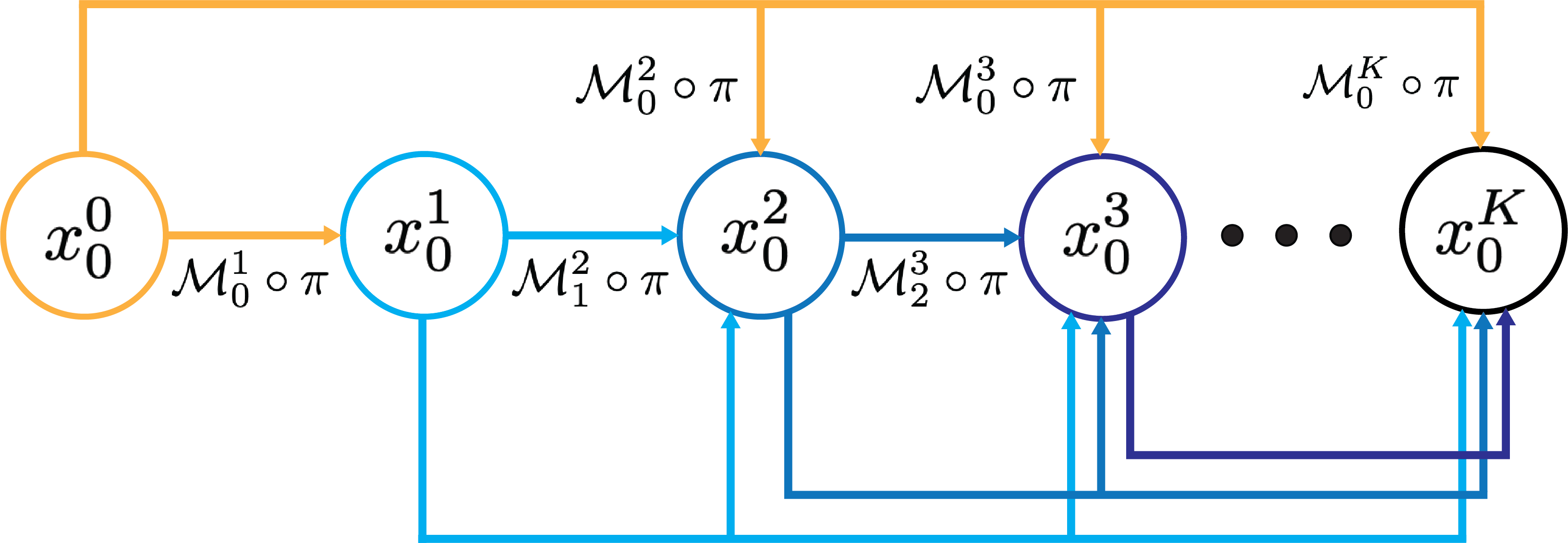} 
            \small
        \end{tikzonimage}
        \caption{Sequential structure of conditional global search paths}
        \label{fig: search tree}
    \end{figure}
}
{
}

As described above, the initial guess for problem $\mathcal{P}_k$ can be generated using solutions from any previously solved simpler problem $\mathcal{P}_{k'}$. 
These different options can be visualized as a network of projection pathways, as illustrated in Fig. \ref{fig: search tree}, where each connection represents a distinct mapping operator. 
The hierarchical structure demonstrates that as $k$ increases, both the number of available source problems $\mathcal{P}_{k'}$ and the variety of possible mapping operators $\mathcal{M}_{k'}^k$ expand, providing greater flexibility in initialization strategy selection.

Additionally, for any map $\mathcal{M}_{k'}^k$, there can be many different ways to generate initial guesses for problem $\mathcal{P}_k$ using the solutions from previously solved problems $\mathcal{P}_{k'}$.
Since the solution space $\mathcal{X}_k$ is higher-dimensional than $\mathcal{X}_{k'}$, the mapping $\mathcal{M}_{k'}^k$ is generally not unique, offering multiple possible strategies for constructing initial guesses in $\mathcal{X}_k$ using solutions in $\mathcal{X}_{k'}$.
Prior to further discussion on the various strategies, it is essential to recognize that each optimal solution $x^{*, k'}$ to $\mathcal{P}_{k'}$ comprises both reference and realization components.
To construct the initial guess for the reference solution for $\mathcal{P}_k$, a natural strategy might be to set the reference solution for $\mathcal{P}_{k'}$ identical to the reference solution for $\mathcal{P}_{k}$, which is what is done in the subsequent numerical experiments in this study.
However, the initial guess generation protocol for the realization solutions in $\mathcal{P}k$ admits multiple options.
A simple approach is to assign the reference solution from $\mathcal{P}_{k'}$ identically to each realization solution in $\mathcal{P}_k$, resulting in a single configuration.
An alternative strategy could be to map the realization solutions for $\mathcal{P}_{k'}$ to the realization solutions for $\mathcal{P}_{k}$.
However, this latter approach introduces additional complexity due to the dimensional mismatch between $\mathcal{X}_{k'}$ and $\mathcal{X}_{k}$.
We provide a comprehensive discussion on this approach in the following paragraph.

Let $x^{*, k'}$ denote a locally optimal solution from previously solved problem $\mathcal{P}_{k'}$ and $x_0^{k}$ denote the initial guess we are trying to construct for $\mathcal{P}_{k}$.
For a given random sample $\omega$, the corresponding realization solution $x^{*, k', \omega}$ constitutes a subset of the complete solution $x^{*, k'}$, denoted as $x^{*, k', \omega} \subset x^{*, k'}$.
Similarly, for random sample $\omega$, the corresponding realization component $x_{0}^{\omega}$ forms a subset of the complete initial guess $x_{0}^{k}$, denoted as $x_{0}^{\omega} \subset x_{0}^{k}$.
The total number of possible mapping operators $\mathscr{C}_{k'}^k$ between solutions $x^{*, k'}$ and $x_0^{k}$ is given by:
\begin{equation}
    \mathscr{C}_{k'}^k = k'^{k} + 1,
\end{equation}
where the first term enumerates all possible realization-to-realization mappings from $\mathcal{P}_{k'}$ to $\mathcal{P}_k$, while the unity term accounts for the singular reference-to-realization mapping.
$\mathscr{C}_{k'}^k$ comprehensively captures all possible assignments, including cases where some realization solutions are not assigned to any initial guess, or where multiple initial guesses are constructed using the same realization solution.
An alternative approach could also involve mapping the realization solutions of $\mathcal{P}_{k'}$ to the reference solution of $\mathcal{P}_{k}$, further increasing the number of possible strategies. 
However, this study restricts consideration to maps that preserve the reference solution structure, i.e., where the reference solution of $\mathcal{P}_{k'}$ maps directly to that of $\mathcal{P}_k$, as this maintains solution coherence across depths of robustness.

To provide further clarity and illustrate the realization-to-realization mapping strategies, the following examples illustrate how $\mathscr{C}_{k'}^k$ varies with different combinations of source ($k'$) and target ($k$) robustness depths.

\noindent
\textbf{Example 1: $k' = 1$, $k = 2$} \\ 
When $x^{*, 1}$ contains only one realization solution $\bm{x}^{*, 1, \omega_0}$, all realizations in $x_0^2$ must be constructed with this realization. 
Thus, there is only one possible combination for $x_0^2$ given by $(\bm{x}^{*, 1, \omega_0}, \bm{x}^{*, 1, \omega_0})$, and the total number of mappings is $\mathscr{C}_{k'}^k = 2$, which includes the additional reference-to-realization mapping.

\noindent
\textbf{Example 2: $k' = 2$, $k = 3$} \\ 
When $x^{*, 2}$ contains two realization solutions $\bm{x}^{*, 2, \omega_0}$ and $\bm{x}^{*, 2, \omega_1}$, each of the three realization initial guesses in $x_0^3$ can be independently constructed using two of these solutions. 
The total number of mappings is $\mathscr{C}_{k'}^k = 9$, a subset of which includes 
$\{
(\bm{x}^{*, 2, \omega_0}, \bm{x}^{*, 2, \omega_0}, \bm{x}^{*, 2, \omega_0}), 
(\bm{x}^{*, 2, \omega_0}, \bm{x}^{*, 2, \omega_0}, \bm{x}^{*, 2, \omega_1}), 
\dots, 
(\bm{x}^{*, 2, \omega_1}, \bm{x}^{*, 2, \omega_1}, \bm{x}^{*, 2, \omega_1})
\}$.

\noindent
\textbf{Example 3: $k' = 3$, $k = 4$} \\ 
When $x^{*, 3}$ contains three realization solutions $\bm{x}^{*, 3, \omega_0}$, $\bm{x}^{*, 3, \omega_1}$ and $\bm{x}^{*, 3, \omega_2}$, each of the four realization initial guesses in $x_0^3$ can now be independently constructed using three of these solutions. 
The total number of mappings is $\mathscr{C}_{k'}^k = 82$, a subset of which includes: 
$\{
(\bm{x}^{*, 3, \omega_0}, \bm{x}^{*, 3, \omega_0}, \bm{x}^{*, 3, \omega_0}, \bm{x}^{*, 3, \omega_0}),
(\bm{x}^{*, 3, \omega_0}, \bm{x}^{*, 3, \omega_0}, \bm{x}^{*, 3, \omega_0}, \bm{x}^{*, 3, \omega_1}),
\dots, \\
(\bm{x}^{*, 3, \omega_2}, \bm{x}^{*, 3, \omega_2}, \bm{x}^{*, 3, \omega_2}, \bm{x}^{*, 3, \omega_2})
\}$.

For practical implementation, we employ the most straightforward realization mapping strategy, which is to map the reference solution for $\mathcal{P}_{k'}$ onto each realization solution for $\mathcal{P}_k$.
Since the number of control segments in the reference solution differs from the number of segments in the realization solution due to the adaptive segmentation strategy (i.e., $N^\dagger \neq N^\omega$), adjustments to the mapping procedure are necessary to ensure compatibility. 
For the thrust vector components, we establish the following protocol.
The reference control segments preceding the missed thrust event ($\bm{u_p^\dagger}$ where $1 \leq p < n_{\omega}$) are excluded from the mapping, where $n_{\omega}$ denotes the reference segment index at which the missed thrust event initiates. 
The remaining control segments ($\bm{u_p^\dagger}$ with $n_{\omega} \leq p \leq N^\dagger$) are then used to construct the initial guess for the realization’s thrust vector parameters, preserving structural similarity between the two solutions. 
Additionally, the initial coast time ($T_i^\omega$) is set to zero during initialization of the local search to provide a better starting guess for the realization solution, which should begin thrusting immediately after diverging from the reference solution. 
It should be noted that while $T_i^\omega$ is initialized to zero, its optimization bounds remain unrestricted, as specified in Table \ref{tab: bounds_nonconditional_global_search}, allowing full exploration of the feasible space during optimization.

A comprehensive analysis of results obtained from these initialization strategies is presented in \S \ref{section: results and discussion}.
\section{Dynamical Model}
\label{section: dynamical model}

During the preliminary mission design phase, leveraging a high-fidelity dynamical model to design LT trajectories ensures a closer alignment with the actual mission dynamics. 
This process requires numerically solving differential equations that account for the gravitational influences of all relevant celestial bodies, while incorporating their ephemeris data to accurately represent their positions and velocities over time.
Generally, the N-body dynamical model describes the motion of the spacecraft $P$ within an inertial frame $\mathscr{I}$ with respect to a central body $\oplus$ under the gravitational effect of the same central body along with other perturbing bodies $Q_i$ given by:
\begin{align}
    \label{eq: n-body dynamics}
    \ddot{\bm{q}}_{P/\oplus} = -G (m_P + m_\oplus) \frac{\bm{q}_{P/\oplus}}{|\bm{q}_{P/\oplus}|^3} +G \sum_{i=0}^{N_Q-1} m_{Q_i}\left(\frac{\bm{q}_{P/Q_i}}{|\bm{q}_{P/Q_i}|^3}-\frac{\bm{q}_{\oplus/Q_i}}{|\bm{q}_{\oplus/Q_i}|^3}\right) + \langle \frac{\bm{u}}{m_P} , \bm{\hat{q}}_{P/\oplus} \rangle
\end{align}
where $\bm{q}_{Q_j/Q_k}$ denotes the position vector of $Q_j$ with respect to $Q_k$, $m_{Q_i}$ is the mass of the celestial body $Q_i$, $N_Q$ represents the total number of perturbing celestial bodies, $u$ is the control input or thrust vector, $\bm{\hat{q}}_{P/\oplus}$ represents a unit vector from the spacecraft $P$ to the central body $\oplus$, and G is the universal gravitational constant.
In the context of this study, the mass $m_\oplus$ represents the Earth; the mass $m_P$ represents the Power and Propulsion Element (discussed in more detail in \S \ref{section: problem setup}), and the additional masses $m_{Q_i}$ correspond to other bodies that exert a non-negligible gravitational attraction.
Since the spacecraft, at various points along its trajectory, makes close approaches to the Moon, it is included in the N-body ephemeris model, along with the Sun and Jupiter (since the solar and jovian gravitational effects have a non-negligible influence on trajectories in the cislunar space).
Although the higher-order terms in Earth's/Moon's gravity field can potentially impact the trajectory, for the sake of simplicity, the model was confined to contain only point-mass contributions from these bodies.
The relative position of each perturbing body with respect to the central body $\bm{q}_{\oplus/Q_i}$ is computed instantaneously using ephemeris data available through the Spacecraft, Planet, Instrument, C-matrix, and Events (SPICE) database developed by the Navigation and Ancillary Information Facility (NAIF) at NASA \cite{acton_1996_ancillary}.

For a LT trajectory, it is also necessary to account for the change in the spacecraft mass, which can be done by simply augmenting the mass to the state of the spacecraft, where the change in the mass $m$ is governed by the differential equation:
\begin{align} 
\label{eq: mass dynamics}
    \dot{m} = - \frac{|\bm{u}|}{I_{sp} g}
\end{align}
where $|\bm{u}|$ is the 2-norm and hence the thrust magnitude, $g = 9.806 \, \text{m/s}^2$ is the gravitational acceleration on Earth and $I_{sp}$ is the constant specific impulse of the propulsion system.
We neglect other perturbations on the spacecraft e.g., solar radiation pressure, such that the only other term affecting the dynamics is the effect of the control input.
\section{Problem Setup}
\label{section: problem setup}

In this section, we present the Power and Propulsion Element as an LT case study to evaluate and compare the efficacy of different initial guess strategies within the context of a realistic global search mission scenario.
As part of NASA’s overarching objective to extend human exploration beyond the Low Earth Orbits and ultimately reach Mars, the agency is developing the Lunar Orbital Platform-Gateway (or simply, Lunar Gateway/LG), a core component of the Artemis program. 
Envisioned as a modular outpost assembled in lunar orbit, the LG will facilitate travel to and from the lunar surface, and support critical interplanetary missions. 
Central to this effort is the Power and Propulsion Element (PPE), a solar electric propulsion module that provides station-keeping, orbital transfers, and space-tug capabilities to the LG. 
Scheduled to launch alongside the Habitation and Logistics Outpost in late 2025, its success is not only crucial for the Artemis program but also for the broader goals of establishing a sustainable human presence in cislunar space and beyond.

Designing long-duration LT transfers within multibody cislunar environments poses substantial challenges, particularly when MTEs disrupt critical maneuvers. 
The current baseline transfer for the PPE exhibits multiple long thrust arcs, one which lasts approximately 276 days (encompassing both the \emph{spiral} and \emph{alignment} phases; see \S \ref{subsection: problem setup: mission phases} for details on these phases) \cite{mcguire_overview_2021}.
By extrapolating empirical probabilities of safe mode events from previous low-thrust mission data, we estimate a 33\% likelihood of experiencing a single MTE during this thrust arc.
An MTE can potentially result in significant deviations from the nominal trajectory, jeopardizing its ability to reach the target operational orbit, or result in significant increases in fuel consumption and/or flight time (as shown in the motivating example in Fig. \ref{fig:motivating_example}). 
Such an event can severely impact the mission by increasing fuel consumption or extending the overall flight time, warranting a more comprehensive understanding of the robust LT solution space, and therefore more efficient global search algorithms for robust LT trajectory design.

In recent years, significant attention has been devoted to transfers toward key cislunar orbits, with particular emphasis on the NRHO. Prior studies have examined impulsive \cite{pritchett_impulsive_2018, whitley_earth-moon_2018} and LT \cite{mcguire_lowthrust_2017, pritchett_impulsive_2018} transfers from (super) geostationary orbits, as well as ballistic trajectories leveraging multibody dynamical structures \cite{parker_low_2007, parker_modeling_2008, parker_monthly_2010, parker_survey_2011}. Recent efforts have focused on eclipse-conscious trajectories \cite{singh_eclipse_2021, pascarella_low_2023} and station-keeping analyses \cite{davis_orbit_2017, guzzetti_stationkeeping_2017}, while there is also growing interest in transfers between NRHOs and other cislunar orbits \cite{mccarty_missed_2020, oshima_use_2019, zhang_transfers_2020}, including lunar landing strategies \cite{trofmov_near-rectilinear_2020}.
A recent study by Karn et al. focuses on understanding the impact of an MTE on the LG, but their focus is exclusively on the terminal insertion phase \cite{karn_recovery_2024}.
Despite recent advances, the global robust solution landscape topology for LT cislunar transfers remains relatively underexplored, making effective initial guess generation challenging.

\subsection{Orbital Parameters}
\label{subsection: problem setup: orbital parameters}

The target operational orbit chosen for the LG is the Earth-Moon $\mathcal{L}_2$ Southern Near Rectilinear Halo Orbit with a 9:2 lunar synodic resonance \cite{whitley_options_2015, lee_gateway_2019}. 
For simplicity, we will simply refer to this orbit as the \textit{NRHO} (i.e., in singular form) throughout the rest of the paper.
The NRHO has an average perilune radius of 3,366 km, with a minimum altitude of 1,450 km over the northern lunar hemisphere. 
The apolune radius extends to about 70,000 km, with a minimum altitude of 68,000 km over the southern lunar hemisphere. 
The orbital period is approximately 6.56 days.
It provides distinct advantages such as low orbit maintenance costs, convenient access to other key cislunar orbits as well as the lunar surface, and minimal eclipse durations.
It is sufficiently stable, with an annual station-keeping $\Delta v$ budget of approximately 15 mm/s \cite{davis_orbit_2017}. 
These characteristics make it an ideal choice for LG's operational orbit.

\begin{table}[!htbp]
\begin{minipage}[b]{0.45\textwidth}
\centering
\caption{Geostationary Transfer Orbit}
\label{tab: sGTO}
\begin{tabular}{ll}
\hline
Parameter                     & Value          \\ \hline
Period {[}days{]}             & 26.41          \\
Mean Periapse Radius {[}km{]} & 6,578          \\
Mean Apoapse Radius {[}km{]}  & 40,278        
\end{tabular}
\end{minipage}
\hfill
\begin{minipage}[b]{0.5\textwidth}
\centering
\caption{Reference $\mathcal{L}_2$ 9:2 Southern NRHO}
\label{tab: NRHO}
\begin{tabular}{ll}
\hline
Parameter                     & Value          \\ \hline
Period {[}days{]}             & 6.56           \\
Mean Perilune Radius {[}km{]} & 3,366          \\
Mean Apolune Radius {[}km{]}  & 71,000        
\end{tabular}
\end{minipage}
\end{table}

The baseline trajectory for the transfer of the PPE module, referred to as the Design Reference Mission (DRM), has undergone several iterations. 
Our problem formulation aligns closely with the third iteration, DRM-3, which begins from a Geostationary Transfer Orbit (GTO) and terminates at the NRHO. 
The initial orbit is a 26.41-day orbit characterized by a mean periapse radius of 6,578 km (periapse altitude of 207 km) and a mean apoapse radius of 40,278 km (apoapse altitude of 33,907 km) (Table~\ref{tab: sGTO}). 
The terminal orbit is a 6.56-day NRHO with a mean perilune radius of 3,366 km (periapse altitude of 207 km) and a mean apolune radius of 71,000 km (apoapse altitude of 68,000 km) (Table~\ref{tab: NRHO}).

\begin{table}[!htbp]
\centering
\caption{Spacecraft Parameters}
\label{tab:spacecraft}
\begin{tabular}{ll}
\hline
Parameter                         & Value                                       \\ \hline
Wet Mass {[}kg{]}                 & 15,000                                      \\
Dry Mass {[}kg{]}                 & 5,000                                       \\
Fuel Mass {[}kg{]}                & 10,000                                      \\
Thrust Acceleration {[}$m/s^2${]} & 2 $\times$ 10\textsuperscript{-4}           \\
Specific Impulse {[}s{]}          & 2,708         
\end{tabular}
\end{table}

The spacecraft parameters were selected to closely align with the characteristic properties of the PPE (Table \ref{tab:spacecraft}).
The stack wet mass is set at 15,000 kg, comprising a dry mass of 5,000 kg and a fuel mass of 10,000 kg. 
The propulsion system, which includes a combination of four Busek BHT-6000 thrusters and two NASA AEPS thrusters, provides a maximum thrust acceleration of 2 $\times$ 10\textsuperscript{-4} m/s\textsuperscript{2} and a constant specific impulse of 2,708 s, with the maximum thrust acceleration defined as the ratio of the maximum attainable thrust to the initial wet mass.

\subsection{Mission Phases}
\label{subsection: problem setup: mission phases}

The DRM-3 framework consists of four stages, each playing a critical role in the mission's design \cite{mcguire_power_2020, mcguire_overview_2021}. 
In the spiral phase, the spacecraft employs a continuous thrust, directed anti-parallel to its instantaneous velocity vector, to gradually raise its altitude. 
This is followed by the alignment phase, during which the thrust vector is optimized to achieve a fuel-efficient trajectory, ensuring the spacecraft is properly oriented for the next phase. 
Next, in the ballistic phase, the spacecraft coasts along a precisely designed ballistic arc toward the Moon. 
Finally, the mission culminates in the insertion phase, which involves a precise thrust maneuver to inject the spacecraft into the NRHO. 
In this study, we adopt a similar conceptual framework but simplify the trajectory design by consolidating these stages into two primary phases.

\subsubsection{Spiral Phase}
\label{subsubsection: problem setup: mission phases: spiral phase}

The \emph{spiral phase} in our framework closely resembles that of DRM-3, with the primary difference being the initial epoch, which is set to November 1, 2025, to align with the anticipated launch date. 
The spacecraft applies continuous thrust anti-parallel to its instantaneous velocity vector, gradually raising its altitude. 
Over the course of approximately 178 days, the spacecraft spirals outward to a distance of about 25.5 Earth radii, achieving a Jacobi integral value by the end of this phase that matches that of NRHO.
During this process, the spacecraft consumes approximately 1,740 kg of fuel. 
For simplicity, we exclude the occurrence of any MTEs during this phase.

\subsubsection{Low-Thrust Transfer Phase}
\label{subsubsection: problem setup: mission phases: LT transfer phase}

The \emph{LT transfer phase} consolidates the alignment, ballistic, and insertion phases into a single optimization problem. 
This phase is formulated as a minimum-fuel optimal control problem, where the initial boundary condition is the terminal state of the spiral phase, and the terminal boundary condition is the NRHO. 
The optimizer is permitted to select the insertion point along the NRHO, introducing an additional degree of freedom to correct for any necessary phasing adjustments.
To facilitate this optimization, DyLAN is used to supply a binary space partitioning (\texttt{bsp}) file containing the ephemeris data for the target orbit. 
The ephemeris is interpolated using splines, enabling the computation of gradients for the objective function with respect to the insertion epoch. 
This allows the optimizer to adjust the insertion epoch dynamically, identifying the optimal insertion point that minimizes fuel consumption.

\subsubsection{Numerical Implementation Details}
\label{subsubsection: problem setup: mission phases: numerical implementation details}

In the non-robust case ($\mathcal{P}_{0}$), the transfer trajectory is discretized into 50 segments using a forward-backward shooting control transcription.
The number of segments was carefully chosen to ensure adequate control authority, resulting in throttle profiles that exhibit the characteristic bang-on, bang-off structure commonly observed in optimal low-thrust solutions.
The maximum allowable time of flight for the transfer phase is constrained to 150 days, comprising a maximum initial coast time ($T_\textrm{i, max}$) of 30 days, a maximum final coast time ($T_\textrm{f, max}$) of 30 days, and a maximum shooting duration ($T_\textrm{s, max}$) of 90 days. 
For the robust case ($\mathcal{P}_{k>0}$), the reference and realization solutions are parameterized differently. 
Notably, while the reference control solution utilizes the same number of segments and decision variable bounds as in the non-robust case, the number of segments and decision variable bounds for the realization solutions are adjusted according to the procedure in \S \ref{subsection: problem formulation: transcription into a nonlinear program}.
The thrust vectors for all solutions are represented in spherical coordinates, with their direction parameterized by in-plane and out-of-plane angles, and their magnitude expressed as the throttle, which varies within the range [0, 1].

Both the non-conditional and conditional approaches involve solving a local optimization problem, which necessitates an initial guess to initialize the local search. 
The non-conditional methodology samples initial points from a prescribed, fixed probability distribution with global support over the solution domain.
Specifically, this study implements the non-conditional method $\mathcal{S}(k)$ using a uniform distribution coupled with monotonic basin hopping \cite{wales_global_1997, leary_global_2000}. 
In contrast, the conditional method $\mathcal{S}(k | k')$ constructs initial guesses for problem $\mathcal{P}_{k}$ using solutions from simpler robust problems $\mathcal{P}_{k'}$ obtained via $\mathcal{S}(k')$.
For both strategies, each candidate solution is subsequently refined through the use of SNOPT \cite{gill_snopt_2005}, which is configured to exploit the sparsity of the Jacobian matrix (Fig. \ref{fig: analytic derivatives sparsity pattern}). 
To accommodate the dimensional scaling with depth of robustness $k$, the maximum allowable runtime for SNOPT is scaled by factor $(1 + k)$. 
This scaling ensures adequate convergence time for meaningful strategy comparisons. 
Similarly, the monotonic basin hopping algorithm in $\mathcal{S}(k)$ implements proportional runtime scaling to maintain consistent exploration depth across increasing robustness levels.


We impose the same numerical tolerances for both the non-robust and the robust solutions. 
Matchpoint defect constraint violations are allowed up to $1 \ \text{km}$ for position, $0.1 \ \text{km/s}$ for velocity, and $1 \ \text{kg}$ for mass. 
All numerical experiments were executed on a dedicated compute node with 2.8 GHz Intel Cascade Lake processors, employing parallel computing resources to ensure uniform computational conditions across all test cases.

\subsection{Example Solutions}
\label{subsection: problem setup: Example Solutions}

\ifthenelse{\boolean{includefigures}}
{
    \begin{figure}[!htb]
        \centering
        \begin{subfigure}[t]{0.45\textwidth}
            \centering
            \includegraphics[keepaspectratio, width=\textwidth]{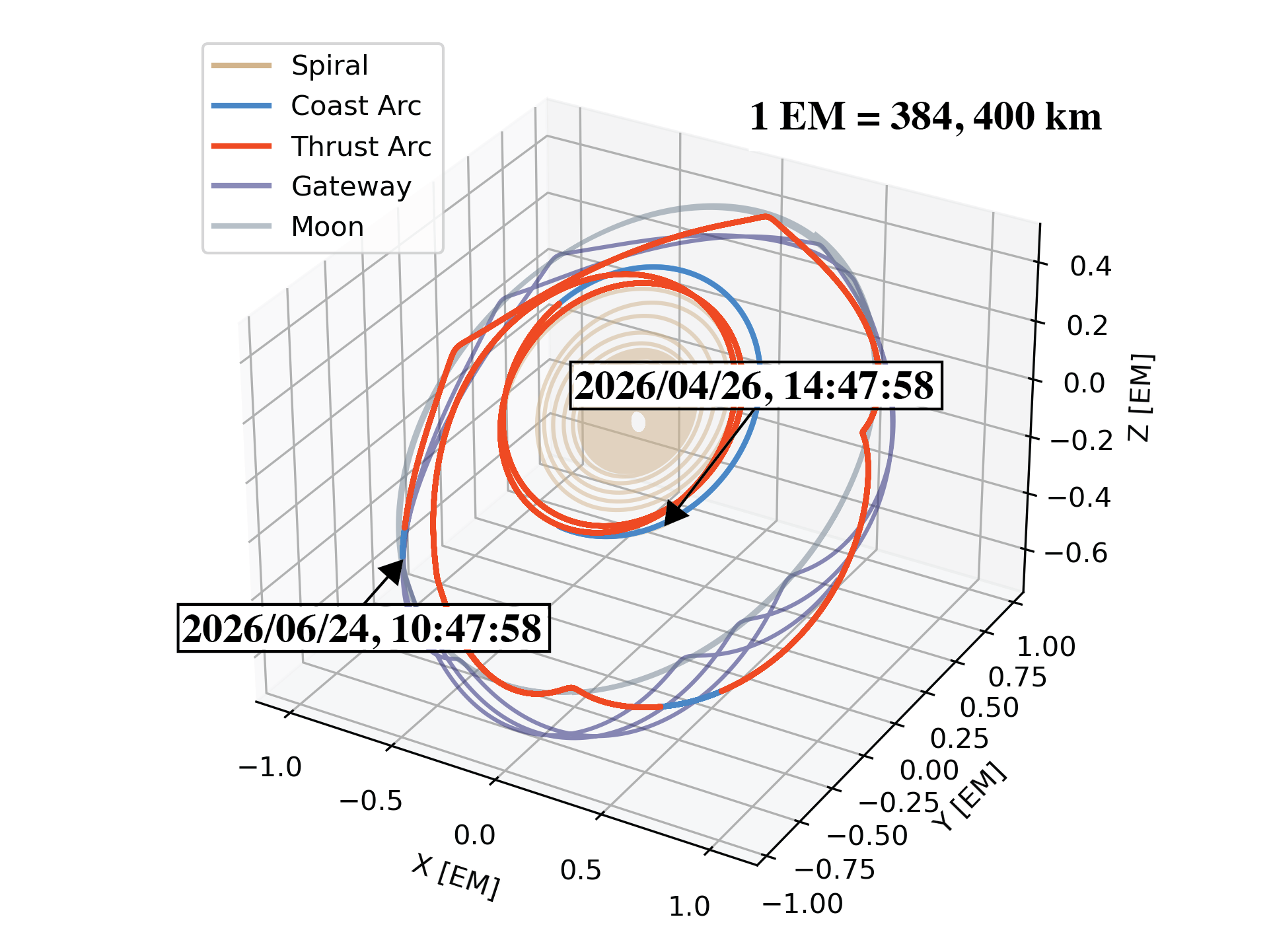}
            \caption{An example non-robust solution}
            \label{fig: non-robust-trajectory}
        \end{subfigure}
        \hfill
        \begin{subfigure}[t]{0.45\textwidth}
            \centering
            \includegraphics[keepaspectratio, width=\textwidth]{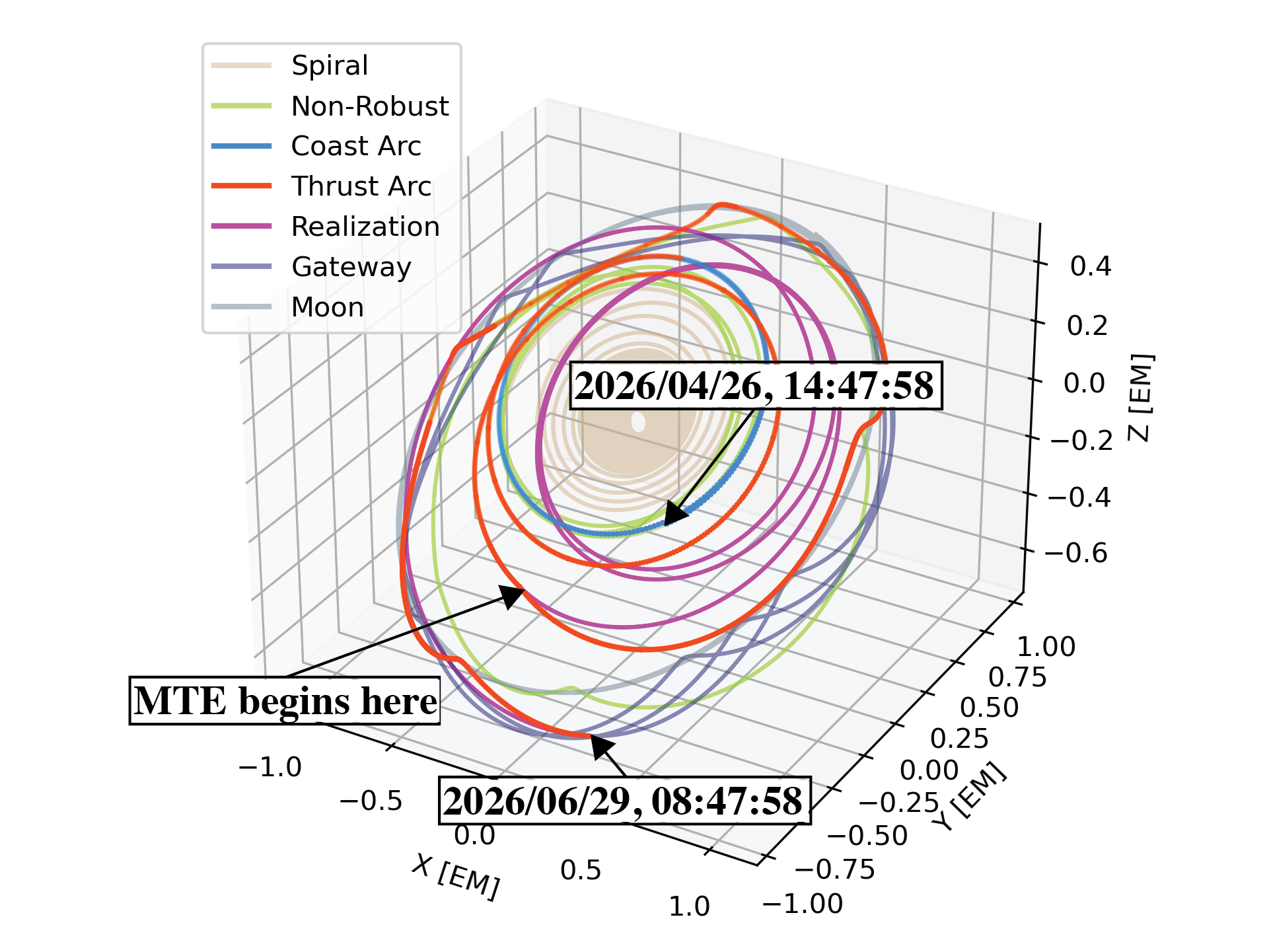}
            \caption{An example robust solution with warm-start}
            \label{fig: robust-trajectory}
        \end{subfigure}
        \caption{Impact of missed thrust events: non-robust vs. robust trajectory solutions.}
        \label{fig: example-solutions}
    \end{figure}
}
{
}

A comparison between an example non-robust and robust trajectories (with one realization), the latter obtained through the conditional search method, is presented in Figure \ref{fig: example-solutions}, demonstrating their fundamental structural differences.
The thrust arcs for the trajectories are shown in \textcolor{Red}{\textbf{red}}, while the coast arcs are shown in \textcolor{DodgerBlue}{\textbf{light blue}}. 
The non-robust solution, shown in Figure \ref{fig: non-robust-trajectory}, serves as initial guess for the robust solution obtained via $\mathcal{S}(1 | 0)$. 
The resulting robust solution, presented in Figure \ref{fig: robust-trajectory}, corresponds to a 60-hour MTE initiated at approximately mid-transfer. 
The robust solution comprises a reference trajectory (thrust arcs in \textcolor{Red}{\textbf{red}}, coast arcs in \textcolor{DodgerBlue}{\textbf{light blue}}) and its realization trajectory (\textcolor{Magenta}{\textbf{magenta}}), with the original non-robust solution shown in \textcolor{Chartreuse}{\textbf{green}} for comparison.
Although the reference trajectory and its realization both insert into the NRHO at the same point near the apolune, they follow significantly different paths to achieve this insertion. 
Making the solution robust introduces a 5-day delay relative to the non-robust trajectory and incurs an additional 35-kg fuel.
\section{Results and Discussion}
\label{section: results and discussion}

The comparative analysis of initialization strategies, namely, \textit{non-conditional global search} and \textit{conditional global search}, employs three fundamental algorithmic performance metrics: \emph{feasibility ratio}, \emph{solving time per solution}, and \emph{solution quality}, which collectively characterize the strategy. 
The feasibility ratio is defined as the percentage of initial guesses that converge to feasible solutions, the solving time per solution is defined as the average compute time required to generate each solution, and the solution quality is quantified by the resulting $\Delta v$ with lower values indicating more fuel-efficient trajectories.
For the non-conditional case, these metrics are computed across the aggregate set of initial guesses generated through basin hopping iterations. 
For the conditional case, these metrics are calculated based on the aggregate set of initial guesses constructed using pertinent prior solutions obtained via the non-conditional approach.
The results highlight the trade-offs between the global search space exploration through the non-conditional approach and the narrower search space exploration through the conditional approach. 
The subsequent sections provide a comprehensive analysis of these strategies, with an emphasis on trends observed in the solutions across different uncertainty realizations.

\begin{table}[!htbp]
\centering
\caption{Numerical Experimental Parameters ($\mathcal{P}_{0}$ and $\mathcal{P}_{1}$)}
\begin{tabular}{lcc}
\hline
                                                                                   & $\mathcal{P}_{0}$         & $\mathcal{P}_{1}$                                                                                     \\ \hline
$\Delta\tau$ {[}days{]}                                                            & 0               & \{0.5\}, \{1.0\}, \{1.5\}                                                                   \\
\begin{tabular}[c]{@{}l@{}}Indices of Segments \\ with MTE Initiation\end{tabular} & N/A             & \{4\}, \{8\}, \{12\}, \{16\}, \{20\}, \{24\}, \{26\}, \{32\}, \{40\}, \{48\}                \\ \hline
\end{tabular}
\label{tab: robust_parameters_k0_k1}
\end{table}

\begin{table}[!htbp]
\centering
\caption{Numerical Experimental Parameters ($\mathcal{P}_{2}$ and $\mathcal{P}_{3}$)}
\begin{tabular}{lcc}
\hline
                                                                                   & $\mathcal{P}_{2}$                                             & $\mathcal{P}_{3}$                                             \\ \hline
$\Delta\tau$ {[}days{]}                                                            & \{0.5\}, \{1.0\}, \{1.5\}                          & \{0.5\}, \{1.0\}, \{1.5\}                          \\
\begin{tabular}[c]{@{}l@{}}Indices of Segments \\ with MTE Initiation\end{tabular} & \begin{tabular}[c]{@{}l@{}} \ \ \{4, 48\}, \\ \{12, 40\}, \\ \{20, 32\}, \\ \ \ \{4, 12\}, \\ \ \ \{8, 16\}, \\ \{12, 20\}, \end{tabular} & \begin{tabular}[c]{@{}l@{}} \ \ \{4, 24, 48\}, \\ \ \ \{4, 26, 48\}, \\ \{12, 24, 40\}, \\ \{12, 26, 40\}, \\ \ \ \ \ \{4, 8, 12\}, \\ \ \ \{4, 12, 20\} \end{tabular} \\ \hline
\end{tabular}
\label{tab: robust_parameters_k2_k3}
\end{table}

Tables \ref{tab: robust_parameters_k0_k1} and \ref{tab: robust_parameters_k2_k3} summarize the MTE parameters used in the numerical experiments for varying robustness complexity, where $k$ denotes the depth of robustness and $\Delta \tau$ specifies MTE duration.
In Table \ref{tab: robust_parameters_k0_k1}, the $\mathcal{P}_{0}$ includes no MTEs, serving as a baseline for comparison. 
For $\mathcal{P}_{1}$, we assume a single realization whose possible initiation locations are chosen to ensure a uniform distribution along the reference trajectory. 
The corresponding outage duration $\Delta\tau$ can be either $0.5$, $1.0$, or $1.5$ days.
For problems with higher depth of robustness $\mathcal{P}_{k>1}$ (Table \ref{tab: robust_parameters_k2_k3}), the duration of each outage is assumed to be identical across all realizations along the trajectory (i.e., we assume each realization undergoes a thruster outage of the same duration), and the segment indices are strategically selected to capture a diverse range of scenarios, with particular emphasis on critical locations along the trajectory, such as the initial, mid-point, and terminal phases of the transfer.
Previous studies examining the relationship between robust solutions and invariant manifolds in three-body problems have demonstrated a stronger sensitivity of low-thrust trajectories during the initial and terminal phases of the transfer \cite{sinha_statistical_2024_arxiv}. 
Although these findings pertain to a different problem, they offer provide meaningful insights for our present investigation.
Building on this insight, we focus the simulation of MTE initiation points toward the beginning and end of the transfer to account for these sensitivities. 
The analysis also incorporates mid-transfer regions corresponding to the $\mathcal{L}_1$ gateway, where solution sensitivity is expected to be amplified by chaotic dynamics. 
This strategic distribution of initiation points ensures comprehensive coverage of dynamically critical regions.

\subsection{Non-Conditional Global Search $\mathcal{S}(k)$}
\label{subsection: results and discussion: non-conditional global search}

\begin{table}[H]
\centering
\caption{Non-Conditional Global Search Performance Metrics}
\begin{tabular}{lcccccccccc}
\hline
                           & $\mathcal{S}(0)$  & \multicolumn{3}{c}{$\mathcal{S}(1)$} & \multicolumn{3}{c}{$\mathcal{S}(2)$} & \multicolumn{3}{c}{$\mathcal{S}(3)$} \\ \hline
$\Delta\tau$ {[}days{]}    & 0    & 0.5    & 1.0    & 1.5   & 0.5    & 1.0    & 1.5   & 0.5    & 1.0     & 1.5   \\ \hline
Number Of Solutions        & 233  & 1053   & 998    & 956   & 161    & 128    & 114   & 8      & 3       & 4     \\
Feasibility Ratio {[}\%{]} & 9.71 & 8.78   & 8.32   & 7.97  & 3.35   & 2.67   & 2.38  & 0.22   & 0.08    & 0.11  \\
Time/Solution {[}h{]}      & 5.15 & 11.40  & 12.02  & 12.55 & 44.72  & 56.25  & 63.16 & 900.00 & 2400.00 & 1800  \\ \hline
\end{tabular}
\label{tab: non-conditional metrics}
\end{table}

Table \ref{tab: non-conditional metrics} presents the performance characteristics of $\mathcal{S}(k)$ across varying depths of robustness ($k$) and missed thrust event durations ($\Delta \tau$).
It is important to emphasize that, although the solutions are categorized by $\Delta \tau$, each category encompasses solutions where the realization may initiate at any of the segments specified in Tables \ref{tab: robust_parameters_k0_k1} and \ref{tab: robust_parameters_k2_k3}.
Analysis reveals that increasing depth of robustness correlates with declining feasibility ratios, reflecting the numerical difficulty in navigating a higher-dimensional solution space, despite the adjustments to maximum allowable runtimes.
For example, for $\Delta \tau = 1.0$ days, the feasibility ratio demonstrates monotonic decay: 8.32\% for $\mathcal{S}(1)$, 2.67\% for $\mathcal{S}(2)$, and 0.08\% for $\mathcal{S}(3)$. 
The computational cost exhibits an inverse behavior, with the average solving time increasing from 12 hours for $\mathcal{S}(1)$ to 56 hours for $\mathcal{S}(2)$, and escalating to 2400 hours for $\mathcal{S}(3)$ at $\Delta \tau = 1.0$ days.

\ifthenelse{\boolean{includefigures}}
{
    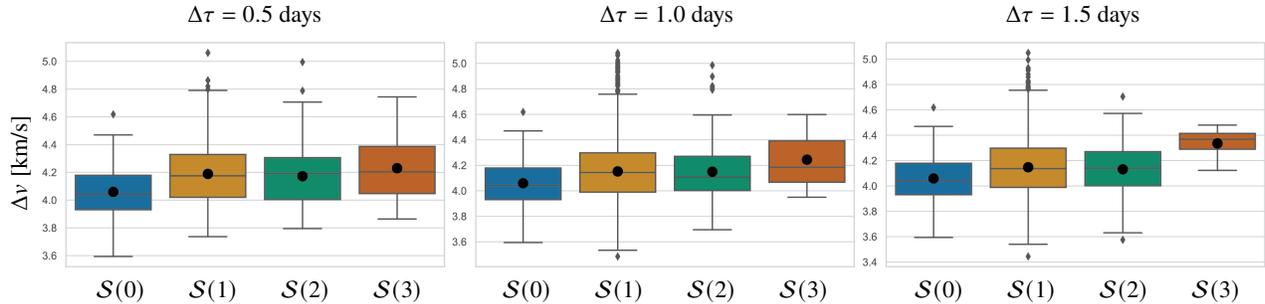
\begin{figure}[!htb]
        \centering
        \begin{tikzonimage}[keepaspectratio, width=0.99\textwidth]{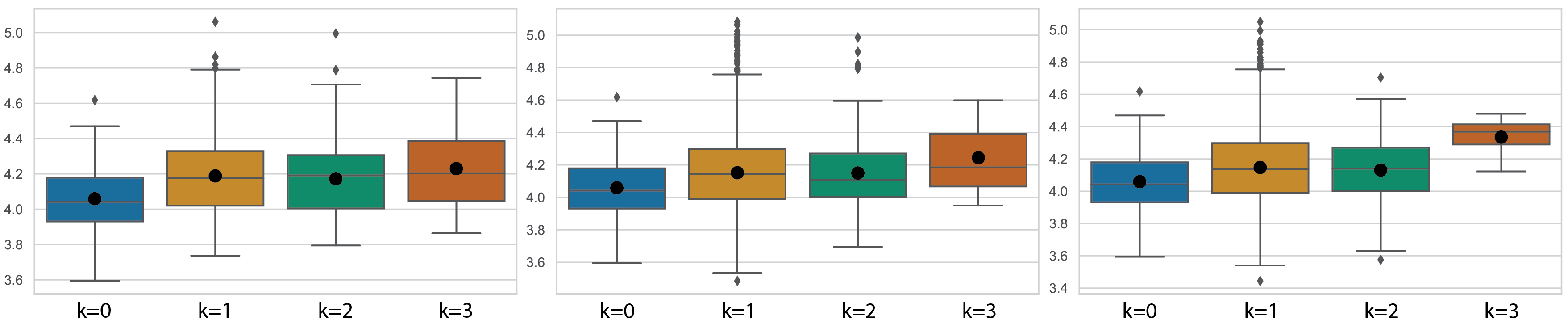}
            \small
            \node[fill=white, opacity=0.0, text opacity=1, anchor=south, rotate=90] at (0.0,0.50) {$\Delta v$ [km/s]};
            \node[fill=white, opacity=0.0, text opacity=1, anchor=south] at (0.175,1.0) {$\Delta \tau = 0.5$ days};
            \node[fill=white, opacity=0.0, text opacity=1, anchor=south] at (0.510,1.0) {$\Delta \tau = 1.0$ days};
            \node[fill=white, opacity=0.0, text opacity=1, anchor=south] at (0.840,1.0) {$\Delta \tau = 1.5$ days};
            \node[fill=white, opacity=1.0, text opacity=1, anchor=north west] at (0.04,0.0725) {$\mathcal{S}(0)$};
            \node[fill=white, opacity=1.0, text opacity=1, anchor=north west] at (0.115,0.0725) {$\mathcal{S}(1)$};
            \node[fill=white, opacity=1.0, text opacity=1, anchor=north west] at (0.19,0.0725) {$\mathcal{S}(2)$};
            \node[fill=white, opacity=1.0, text opacity=1, anchor=north west] at (0.265,0.0725) {$\mathcal{S}(3)$};

            \node[fill=white, opacity=1.0, text opacity=1, anchor=north west] at (0.37,0.0725) {$\mathcal{S}(0)$};
            \node[fill=white, opacity=1.0, text opacity=1, anchor=north west] at (0.445,0.0725) {$\mathcal{S}(1)$};
            \node[fill=white, opacity=1.0, text opacity=1, anchor=north west] at (0.52,0.0725) {$\mathcal{S}(2)$};
            \node[fill=white, opacity=1.0, text opacity=1, anchor=north west] at (0.595,0.0725) {$\mathcal{S}(3)$};

            \node[fill=white, opacity=1.0, text opacity=1, anchor=north west] at (0.705,0.0725) {$\mathcal{S}(0)$};
            \node[fill=white, opacity=1.0, text opacity=1, anchor=north west] at (0.78,0.0725) {$\mathcal{S}(1)$};
            \node[fill=white, opacity=1.0, text opacity=1, anchor=north west] at (0.855,0.0725) {$\mathcal{S}(2)$};
            \node[fill=white, opacity=1.0, text opacity=1, anchor=north west] at (0.93,0.0725) {$\mathcal{S}(3)$};
        \end{tikzonimage}
        \caption{$\Delta v$ (Non-Conditional Global)}
        \label{fig: delta-v non-conditional}
    \end{figure}
}
{
}

Figure \ref{fig: delta-v non-conditional} illustrates the statistical distribution of $\Delta v$ requirements across the same solution categories, where each boxplot highlights the median (horizontal bar), mean (dot), and interquartile range (shaded region), providing quantitative insight into the solution quality.
Robust solutions obtained via $\mathcal{S}(k>0)$ consistently demonstrate higher $\Delta v$ requirements compared to their non-robust counterparts obtained via $\mathcal{S}(0)$, reflecting the additional propellant demands for robust solutions. 
Across $\Delta \tau$, the mean $\Delta v$ requirement exhibits monotonic growth with increasing depth of robustness, with $\mathcal{S}(3)$ showing the most pronounced elevation. 
We observe no strong correlation between $\Delta \tau$ and $\Delta v$, which we hypothesize may result from the relatively narrow range in $\Delta \tau$ values in this study.

\subsection{Conditional global search $\mathcal{S}(k | k')$}
\label{subsection: results and discussion: conditional global search}

\ifthenelse{\boolean{includefigures}}
{
    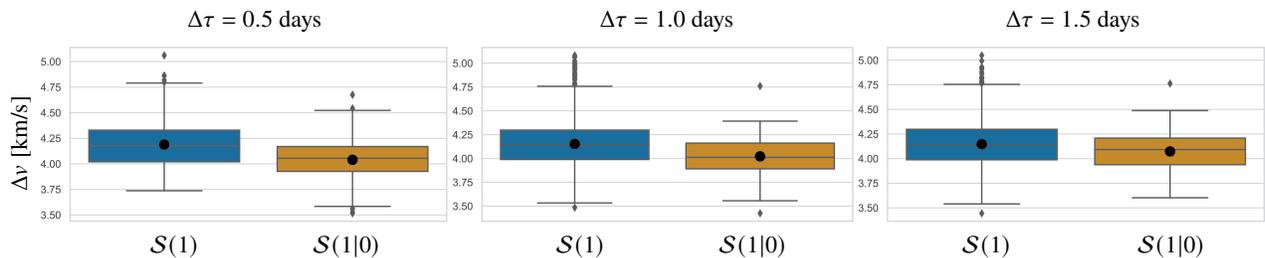
\begin{figure}[!htb]
        \centering
        \begin{tikzonimage}[keepaspectratio, width=0.99\textwidth]{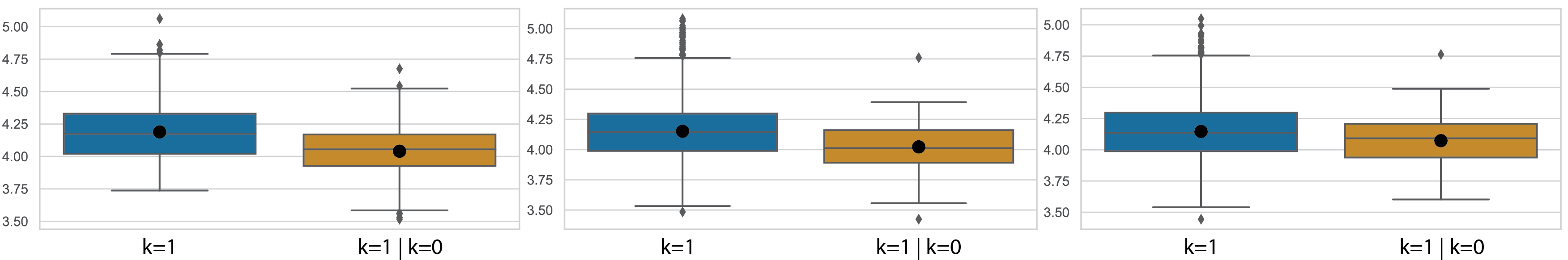}
            \small
            \node[fill=white, opacity=0.0, text opacity=1, anchor=south, rotate=90] at (0.0,0.50) {$\Delta v$ [km/s]};
            \node[fill=white, opacity=0.0, text opacity=1, anchor=south] at (0.175,1.0) {$\Delta \tau = 0.5$ days};
            \node[fill=white, opacity=0.0, text opacity=1, anchor=south] at (0.510,1.0) {$\Delta \tau = 1.0$ days};
            \node[fill=white, opacity=0.0, text opacity=1, anchor=south] at (0.840,1.0) {$\Delta \tau = 1.5$ days};

            \node[fill=white, opacity=1.0, text opacity=1, anchor=center] at (0.11,0.0) {$\mathcal{S}(1)$};
            \node[fill=white, opacity=1.0, text opacity=1, anchor=center] at (0.25,0.0) {$\mathcal{S}(1|0)$};

            \node[fill=white, opacity=1.0, text opacity=1, anchor=center] at (0.44,0.0) {$\mathcal{S}(1)$};
            \node[fill=white, opacity=1.0, text opacity=1, anchor=center] at (0.58,0.0) {$\mathcal{S}(1|0)$};

            \node[fill=white, opacity=1.0, text opacity=1, anchor=center] at (0.77,0.0) {$\mathcal{S}(1)$};
            \node[fill=white, opacity=1.0, text opacity=1, anchor=center] at (0.92,0.0) {$\mathcal{S}(1|0)$};
        \end{tikzonimage}
        \caption{$\Delta v$ (Conditional Global, $\mathcal{P}_1$)}
        \label{fig: delta-v conditional nummte=1}
    \end{figure}
}
{
}

Intuitively, the conditional approach should offer a significant advantage in computational efficiency through a narrower design space constructed by leveraging existing solutions, either non-robust or simpler robust solutions. 
Figure \ref{fig: delta-v conditional nummte=1} presents the total $\Delta v$ distribution for $\mathcal{S}(1 | 0)$, alongside the distribution for $\mathcal{S}(1)$ for comparison. 
$\mathcal{S}(1 | 0)$ demonstrates consistently lower mean $\Delta v$ requirements across all $\Delta \tau$ values, indicating superior solution quality compared to $\mathcal{S}(1)$.
However, the reduced variance in $\mathcal{S}(1 | 0)$ solutions suggests a trade-off between solution quality and diversity, an artifact of constrained domain exploration inherent to the conditional approach. 
Leveraging a diverse collection of prior non-robust solutions could help mitigate this limitation, and improve the diversity of robust solutions obtained via $\mathcal{S}(1 | 0)$.

\ifthenelse{\boolean{includefigures}}
{
    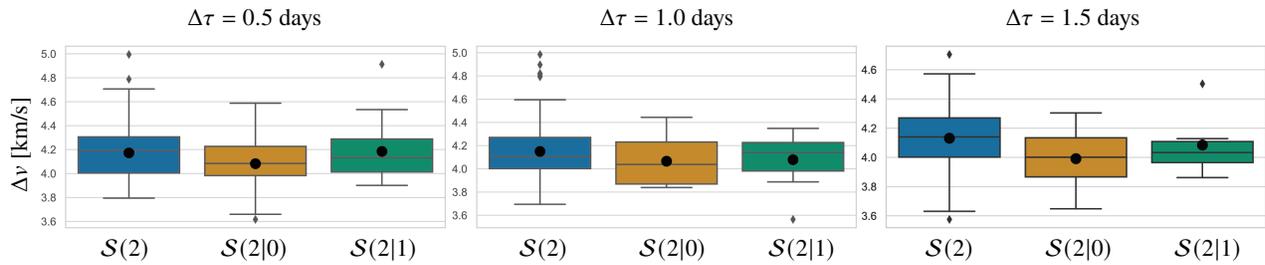
\begin{figure}[!htb]
        \centering
        \begin{tikzonimage}[keepaspectratio, width=0.99\textwidth]{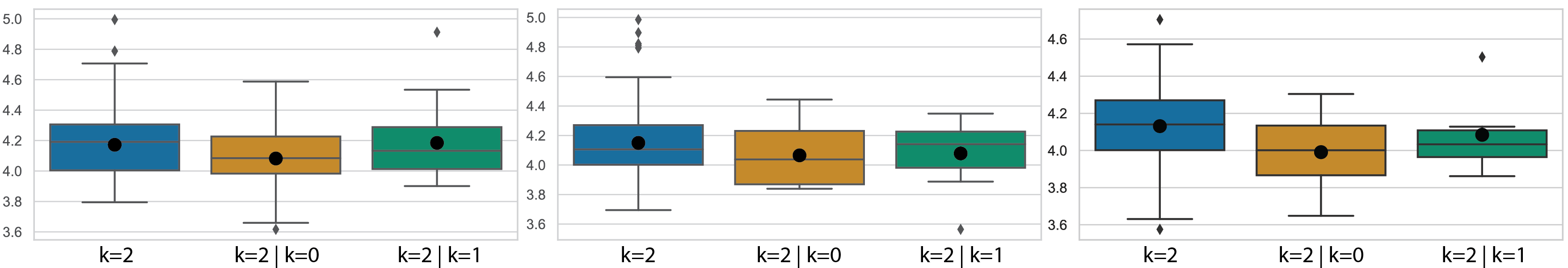}
            \small
            \node[fill=white, opacity=0.0, text opacity=1, anchor=south, rotate=90] at (0.0,0.50) {$\Delta v$ [km/s]};
            \node[fill=white, opacity=0.0, text opacity=1, anchor=south] at (0.175,1.0) {$\Delta \tau = 0.5$ days};
            \node[fill=white, opacity=0.0, text opacity=1, anchor=south] at (0.510,1.0) {$\Delta \tau = 1.0$ days};
            \node[fill=white, opacity=0.0, text opacity=1, anchor=south] at (0.840,1.0) {$\Delta \tau = 1.5$ days};

            \node[fill=white, opacity=1.0, text opacity=1, anchor=center] at (0.07,0.0) {$\mathcal{S}(2)$};
            \node[fill=white, opacity=1.0, text opacity=1, anchor=center] at (0.175,0.0) {$\mathcal{S}(2 | 0)$};
            \node[fill=white, opacity=1.0, text opacity=1, anchor=center] at (0.28,0.0) {$\mathcal{S}(2 | 1)$};

            \node[fill=white, opacity=1.0, text opacity=1, anchor=center] at (0.41,0.0) {$\mathcal{S}(2)$};
            \node[fill=white, opacity=1.0, text opacity=1, anchor=center] at (0.515,0.0) {$\mathcal{S}(2 | 0)$};
            \node[fill=white, opacity=1.0, text opacity=1, anchor=center] at (0.62,0.0) {$\mathcal{S}(2 | 1)$};

            \node[fill=white, opacity=1.0, text opacity=1, anchor=center] at (0.74,0.0) {$\mathcal{S}(2)$};
            \node[fill=white, opacity=1.0, text opacity=1, anchor=center] at (0.845,0.0) {$\mathcal{S}(2 | 0)$};
            \node[fill=white, opacity=1.0, text opacity=1, anchor=center] at (0.95,0.0) {$\mathcal{S}(2 | 1)$};
        \end{tikzonimage}
        \caption{$\Delta v$ (Conditional Global, $\mathcal{P}_2$)}
        \label{fig: delta-v conditional nummte=2}
    \end{figure}
}
{
}

Figure \ref{fig: delta-v conditional nummte=2} compares the $\Delta v$ distributions across three methodologies: non-conditional search $\mathcal{S}(2)$, and two conditional approaches $\mathcal{S}(2 | 0)$ and $\mathcal{S}(2 | 1)$.
Across all $\Delta \tau$ values, $\mathcal{S}(2)$ exhibits consistently higher propellant requirements than the conditional approaches. 
Between the two conditional methods, $\mathcal{S}(2 | 0)$ tends to yield a lower mean $\Delta v$ with a broader range, whereas $\mathcal{S}(2 | 1)$ produces a higher mean $\Delta v$ with a narrower range. 
In this scenario, using the $\mathcal{S}(2 | 0)$ strategy not only improves the solution quality but also leads to more diversity than $\mathcal{S}(2 | 1)$, which can be particularly advantageous in practice during the preliminary mission design phase. 

\ifthenelse{\boolean{includefigures}}
{
    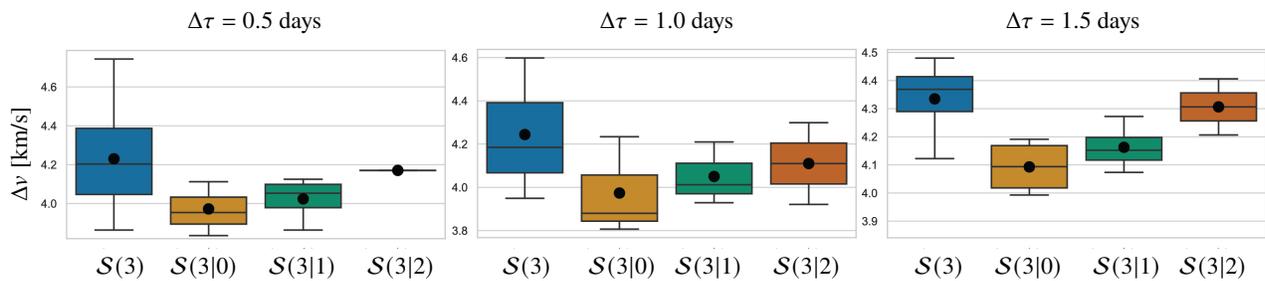
\begin{figure}[!htb]
        \centering
        \begin{tikzonimage}[keepaspectratio, width=0.99\textwidth]{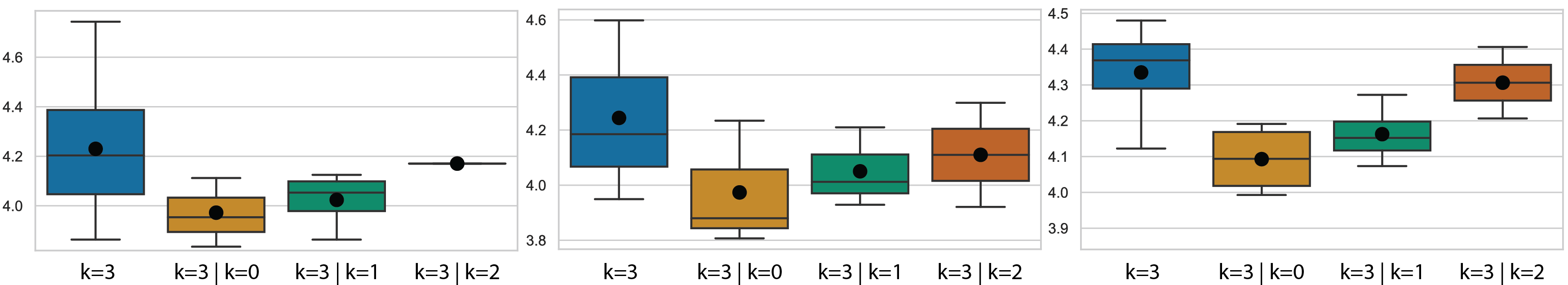}
            \small
            \node[fill=white, opacity=0.0, text opacity=1, anchor=south, rotate=90] at (0.0,0.50) {$\Delta v$ [km/s]};
            \node[fill=white, opacity=0.0, text opacity=1, anchor=south] at (0.175,1.0) {$\Delta \tau = 0.5$ days};
            \node[fill=white, opacity=0.0, text opacity=1, anchor=south] at (0.510,1.0) {$\Delta \tau = 1.0$ days};
            \node[fill=white, opacity=0.0, text opacity=1, anchor=south] at (0.840,1.0) {$\Delta \tau = 1.5$ days};

            \node[fill=white, opacity=1.0, text opacity=1, anchor=north west] at (0.04,0.08) {$\mathcal{S}(3)$};
            \node[fill=white, opacity=1.0, text opacity=1, anchor=north west] at (0.10,0.08) {$\mathcal{S}(3|0)$};
            \node[fill=white, opacity=1.0, text opacity=1, anchor=north west] at (0.18,0.08) {$\mathcal{S}(3|1)$};
            \node[fill=white, opacity=1.0, text opacity=1, anchor=north west] at (0.26,0.08) {$\mathcal{S}(3|2)$};

            \node[fill=white, opacity=1.0, text opacity=1, anchor=north west] at (0.37,0.08) {$\mathcal{S}(3)$};
            \node[fill=white, opacity=1.0, text opacity=1, anchor=north west] at (0.44,0.08) {$\mathcal{S}(3|0)$};
            \node[fill=white, opacity=1.0, text opacity=1, anchor=north west] at (0.51,0.08) {$\mathcal{S}(3|1)$};
            \node[fill=white, opacity=1.0, text opacity=1, anchor=north west] at (0.59,0.08) {$\mathcal{S}(3|2)$};

            \node[fill=white, opacity=1.0, text opacity=1, anchor=north west] at (0.705,0.08) {$\mathcal{S}(3)$};
            \node[fill=white, opacity=1.0, text opacity=1, anchor=north west] at (0.77,0.08) {$\mathcal{S}(3|0)$};
            \node[fill=white, opacity=1.0, text opacity=1, anchor=north west] at (0.85,0.08) {$\mathcal{S}(3|1)$};
            \node[fill=white, opacity=1.0, text opacity=1, anchor=north west] at (0.92,0.08) {$\mathcal{S}(3|2)$};
        \end{tikzonimage}
        \caption{$\Delta v$ (Conditional Global, $\mathcal{P}_3$)}
        \label{fig: delta-v conditional nummte=3}
    \end{figure}
}
{
}


Figure \ref{fig: delta-v conditional nummte=3} compares the $\Delta v$ distributions across four methodologies: non-conditional search $\mathcal{S}(3)$ and three conditional approaches $\mathcal{S}(3 | 0)$, $\mathcal{S}(3 | 1)$, and $\mathcal{S}(3 | 2)$. 
Once again, $\mathcal{S}(3)$ exhibits consistently higher propellant requirements than the conditional approaches across all $\Delta \tau$ values. 
Within the conditional approaches, $\mathcal{S}(3 | 0)$ achieves minimal mean $\Delta v$ requirements.
This result may appear counter-intuitive, as directly seeding initial guesses from non-robust solutions should typically require larger corrections to the solutions, leading to an increasing propellant consumption as $k$ increases.
So, one may expect that solutions to robust problems with higher depths of robustness would serve as more effective initial guesses, as they inherently incorporate greater levels of partial robustness. 
However, in this scenario, $\mathcal{S}(3 | 0)$ demonstrates optimal performance in both solution quality and diversity.
We hypothesize that this phenomenon is an artifact of how the feasible solution space evolves with increasing depth of robustness $k$. 
As $k$ increases, the objective landscape becomes progressively constrained, potentially causing initial guesses from prior robust solutions to map to infeasible regions.
Initial guesses from robust solutions undergo sequential maps through intermediate robustness levels, and these intermediate transformations may lead to infeasible solution mappings.
In contrast, initial guesses derived from non-robust solutions are less constrained and may maintain a higher probability of mapping to feasible regions despite increasing depths of robustness.

These findings elucidate the relationship between the two search approaches and the solution quality in robust problems. 
Among the conditional search methods, using non-robust seeds consistently yield superior results across both solution quality and solution diversity, suggesting that maintaining broader feasible region access outweighs the potential benefits of partial robustness in seed solutions.

\ifthenelse{\boolean{includefigures}}
{
    \begin{figure}[!htb]
        \centering
        \begin{tikzonimage}[keepaspectratio, width=0.99\textwidth]{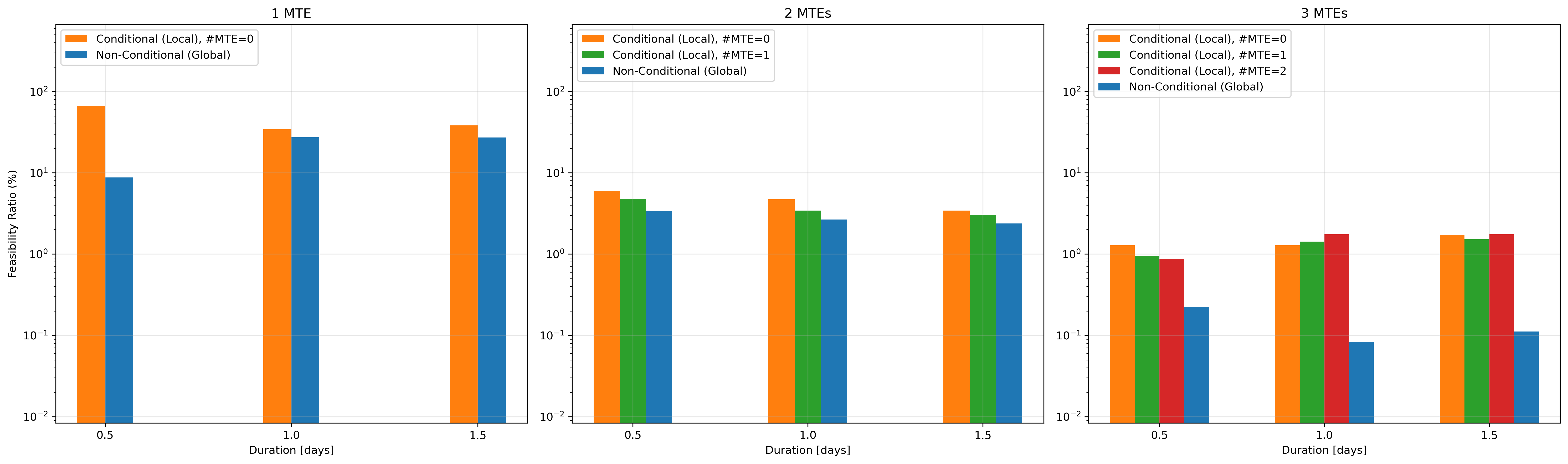}
            \small
            \node[fill=white, opacity=1.0, text opacity=1, anchor=south, rotate=90] at (0.015,0.50) {Feasibility Ratio [\%]};
            \node[fill=white, opacity=1.0, text opacity=1, anchor=south] at (0.19,0.95) {$\mathcal{P}_1$};
            \node[fill=white, opacity=1.0, text opacity=1, anchor=south] at (0.515,0.95) {$\mathcal{P}_2$};
            \node[fill=white, opacity=1.0, text opacity=1, anchor=south] at (0.8475,0.95) {$\mathcal{P}_3$};
            \node[fill=white, opacity=1.0, text opacity=1, anchor=south] at (0.1775,-0.060) {$\Delta\tau$ [days]};
            \node[fill=white, opacity=1.0, text opacity=1, anchor=south] at (0.5125,-0.060) {$\Delta\tau$ [days]};
            \node[fill=white, opacity=1.0, text opacity=1, anchor=south] at (0.8500,-0.060) {$\Delta\tau$ [days]};
            \node[fill=white,opacity=1.0,anchor=center] at (0.1025,0.855){\includegraphics[keepaspectratio, width=0.75 in]{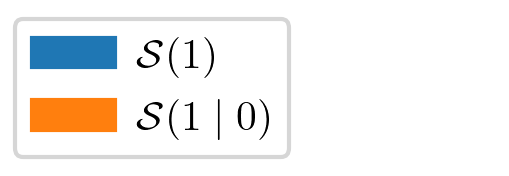}};
            \node[fill=white,opacity=1.0,anchor=center] at (0.43,0.83){\includegraphics[keepaspectratio, width=0.75 in]{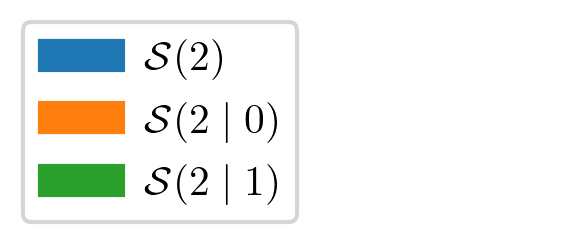}};
            \node[fill=white,opacity=1.0] at (0.7605,0.81){\includegraphics[keepaspectratio, width=0.75 in]{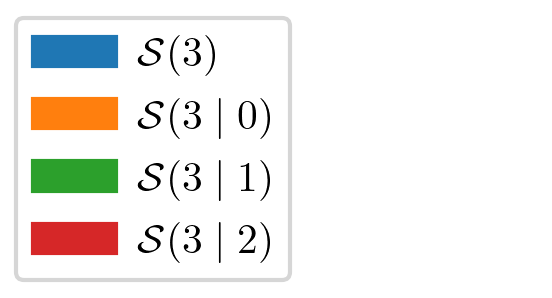}};
        \end{tikzonimage}
        \caption{Feasibility Ratio Comparison for Conditional and Non-Conditional Search Strategies}
        \label{fig: feasibility ratio}
    \end{figure}
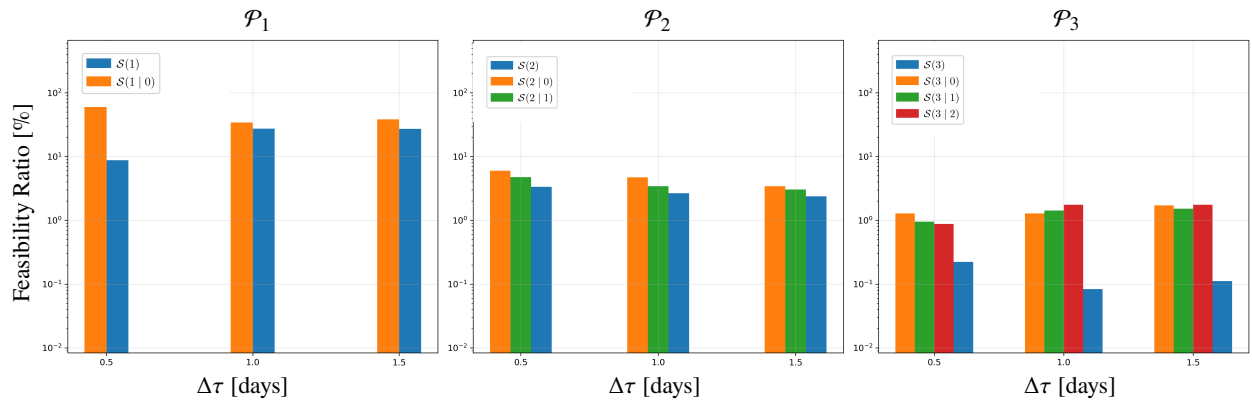
}
{
}

Figure \ref{fig: feasibility ratio} compares the feasibility ratios across three initial-guess generation strategies, namely, a conditional search seeded by non-robust trajectories $\mathcal{S}(k | k'=0)$, a conditional search seeded by robust solutions with lower depths of robustness $\mathcal{S}(k | k')$, and a non-conditional search $\mathcal{S}(k)$.
For $\mathcal{P}_1$, $\mathcal{S}(1 | 0)$ yields a higher feasibility ratios than $\mathcal{S}(1)$ across all $\Delta \tau$ values. 
The transition to $\mathcal{P}_2$ exhibits a significant decrease in feasibility ratios across all methodologies, with conditional approaches maintaining performance advantages over the non-conditional approach. 
Notably, $\mathcal{S}(2 | 0)$ achieves higher feasibility ratios than $\mathcal{S}(2 | 1)$, reinforcing the efficacy of non-robust seeds for problems of increasing depth of robustness.
Finally, for $\mathcal{P}_3$, the feasibility ratios drop even further due to increased dimensionality, yet the conditional search remains more successful than the non-conditional search. 
Overall, from a feasibility perspective, $\mathcal{S}(k | k=0)$ appears to be the optimal initial guess generation strategy.

\ifthenelse{\boolean{includefigures}}
{
    \begin{figure}[!htb]
        \centering
        \begin{tikzonimage}[keepaspectratio, width=0.99\textwidth]{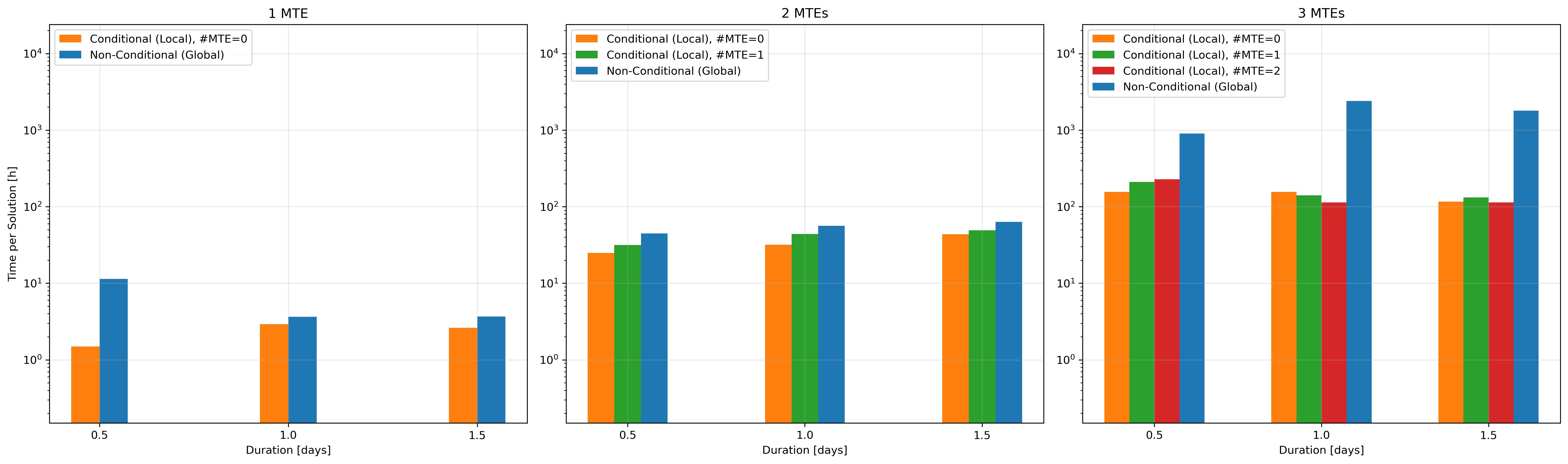}
            \small
            \node[fill=white, opacity=1.0, text opacity=1, anchor=south, rotate=90] at (0.015,0.50) {Time-to-Solve [h]};
            \node[fill=white, opacity=1.0, text opacity=1, anchor=south] at (0.19,0.95) {$\mathcal{P}_1$};
            \node[fill=white, opacity=1.0, text opacity=1, anchor=south] at (0.515,0.95) {$\mathcal{P}_2$};
            \node[fill=white, opacity=1.0, text opacity=1, anchor=south] at (0.845,0.95) {$\mathcal{P}_3$};
            \node[fill=white, opacity=1.0, text opacity=1, anchor=south] at (0.1775,-0.060) {$\Delta\tau$ [days]};
            \node[fill=white, opacity=1.0, text opacity=1, anchor=south] at (0.5125,-0.060) {$\Delta\tau$ [days]};
            \node[fill=white, opacity=1.0, text opacity=1, anchor=south] at (0.8500,-0.060) {$\Delta\tau$ [days]};
            \node[fill=white,opacity=1.0,anchor=center] at (0.10,0.855){\includegraphics[keepaspectratio, width=0.75 in]{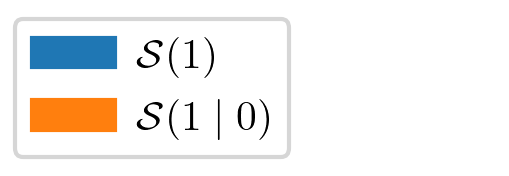}};
            \node[fill=white,opacity=1.0,anchor=center] at (0.43,0.83){\includegraphics[keepaspectratio, width=0.75 in]{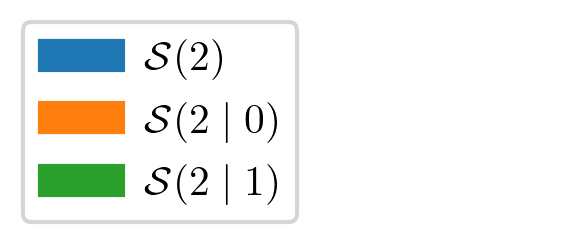}};
            \node[fill=white,opacity=1.0] at (0.7605,0.81){\includegraphics[keepaspectratio, width=0.75 in]{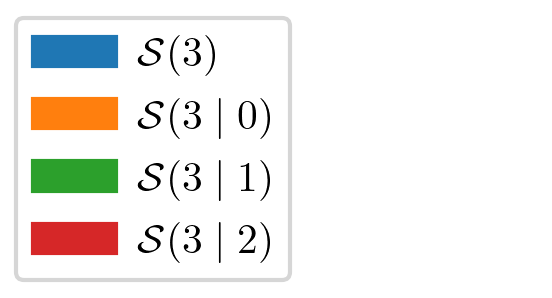}};
        \end{tikzonimage}
        \caption{Average Time-to-Solve Comparison for Conditional and Non-Conditional Search Strategies}
        \label{fig: time to solve}
    \end{figure}
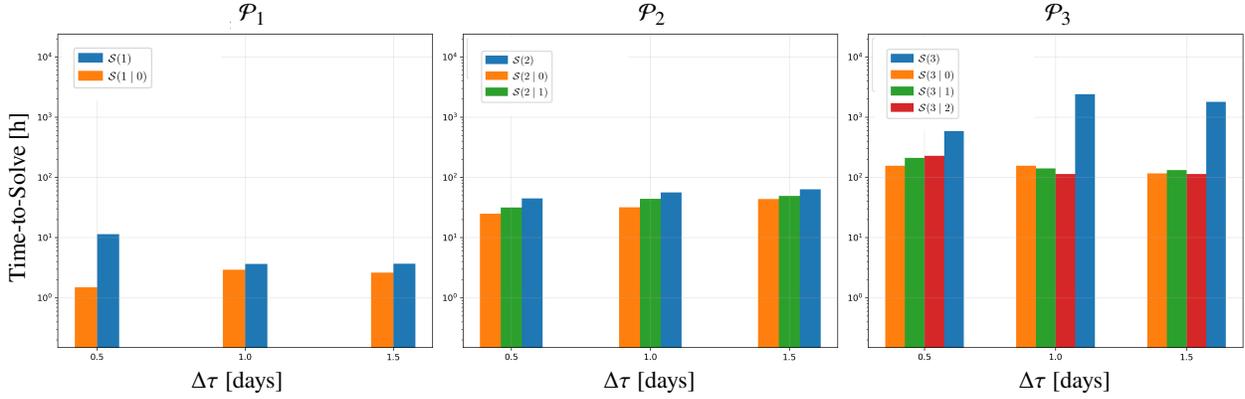
}
{
}

Figure \ref{fig: time to solve} shows the solving time (in hours) for the same three strategies. 
For $\mathcal{P}_1$, $\mathcal{S}(1 | 0)$ appears to be the most efficient, while $\mathcal{S}(1)$ is noticeably slower. 
As we move to $\mathcal{P}_2$, all of the methods require more time due to the higher problem dimensionality, but the conditional approaches, specifically $\mathcal{S}(2 | 0)$, still outperform the non-conditional approach.
Finally, the results for $\mathcal{P}_3$ show further elevation in required solving time, with $\mathcal{S}(3 | 0)$ achieving minimal mean solving times across most $\Delta \tau$ values.
A notable exception occurs at $\Delta \tau = 1.0$ days, where $\mathcal{S}(3 | 2)$ demonstrates marginally better solving times than $\mathcal{S}(3 | 0)$, but this is likely an artifact of the significantly smaller sample size for $\mathcal{S}(2)$ (Table \ref{tab: non-conditional metrics}) and therefore may not be statistically significant.

These results characterize the fundamental trade-offs inherent in the non-conditional approach. 
While enabling comprehensive exploration of the solution space, as evidenced by the large variance in the $\Delta v$ metrics (Figs. \ref{fig: delta-v conditional nummte=1}, \ref{fig: delta-v conditional nummte=2}, \ref{fig: delta-v conditional nummte=3}), the non-conditional approach exhibits diminishing performance with increasing robustness requirements, manifesting in reduced feasibility ratios and escalating solving times.


\ifthenelse{\boolean{includefigures}}
{
    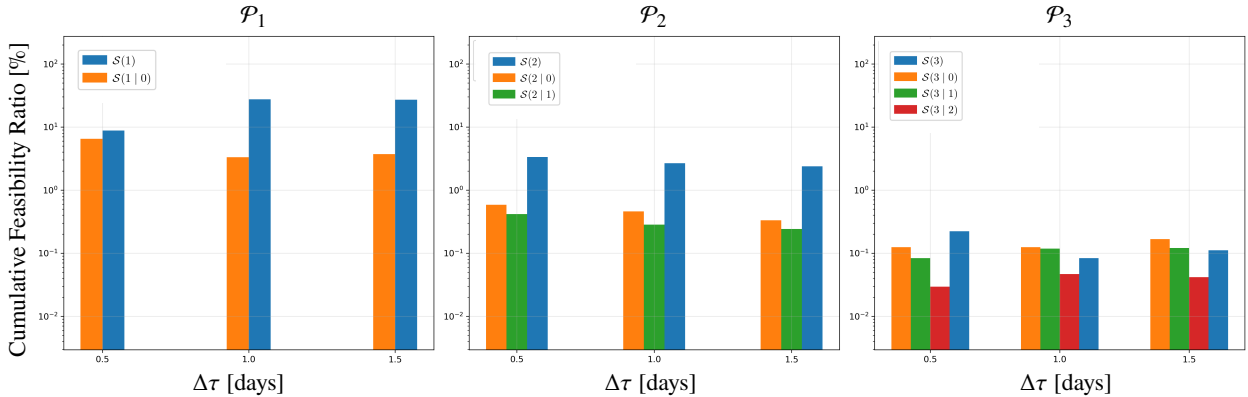
\begin{figure}[!htb]
        \centering
        \begin{tikzonimage}[keepaspectratio, width=0.99\textwidth]{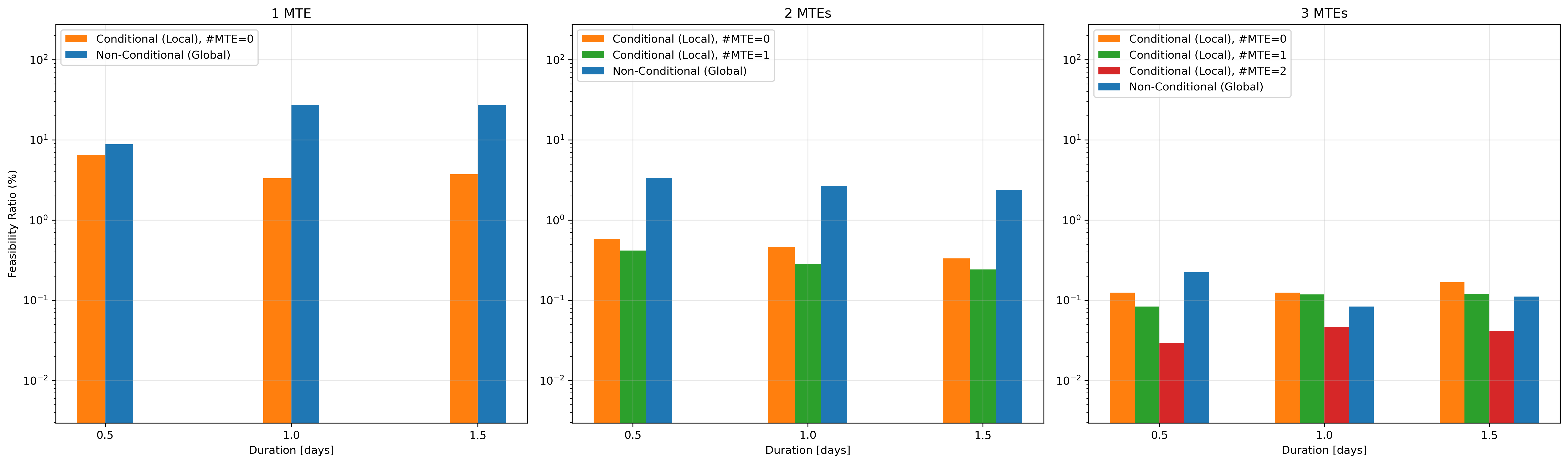}
            \small
            \node[fill=white, opacity=1.0, text opacity=1, anchor=south, rotate=90] at (0.015,0.50) {Cumulative Feasibility Ratio [\%]};
            \node[fill=white, opacity=1.0, text opacity=1, anchor=south] at (0.19,0.95) {$\mathcal{P}_1$};
            \node[fill=white, opacity=1.0, text opacity=1, anchor=south] at (0.515,0.95) {$\mathcal{P}_2$};
            \node[fill=white, opacity=1.0, text opacity=1, anchor=south] at (0.845,0.95) {$\mathcal{P}_3$};
            \node[fill=white, opacity=1.0, text opacity=1, anchor=south] at (0.1775,-0.060) {$\Delta\tau$ [days]};
            \node[fill=white, opacity=1.0, text opacity=1, anchor=south] at (0.5125,-0.060) {$\Delta\tau$ [days]};
            \node[fill=white, opacity=1.0, text opacity=1, anchor=south] at (0.8500,-0.060) {$\Delta\tau$ [days]};
            \node[fill=white,opacity=1.0,anchor=center] at (0.1025,0.855){\includegraphics[keepaspectratio, width=0.75 in]{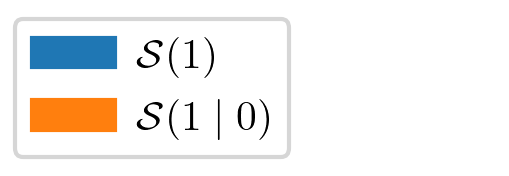}};
            \node[fill=white,opacity=1.0,anchor=center] at (0.435,0.83){\includegraphics[keepaspectratio, width=0.75 in]{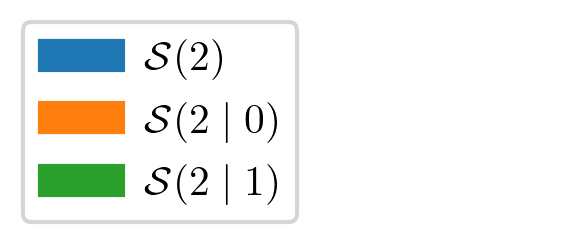}};
            \node[fill=white,opacity=1.0] at (0.7625,0.81){\includegraphics[keepaspectratio, width=0.75 in]{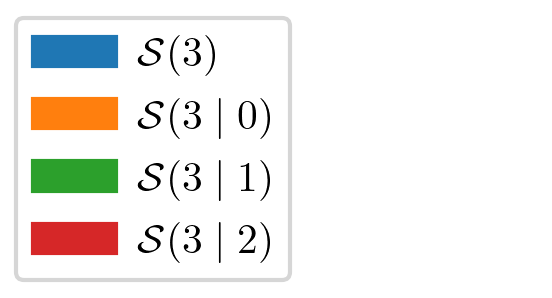}};
        \end{tikzonimage}
        \caption{Cumulative Feasibility Ratio Comparison for Conditional and Non-Conditional Search Strategies}
        \label{fig: adjusted feasibility ratio}
    \end{figure}
}
{
}

When computing the feasibility ratio for the conditional search, it is also essential to incorporate the feasibility ratio of the solutions which informs the initial guesses, necessitating a cumulative metric that combines the feasibility ratios of both the seed generation and conditional search processes.
For instance, if we consider the $\mathcal{S}(1 | 0)$ strategy, the overall feasibility ratio should be the product of the $\mathcal{S}(0)$ feasibility ratio and the $\mathcal{S}(1 | 0)$ feasibility ratio. 
Figure \ref{fig: adjusted feasibility ratio} presents the \emph{cumulative} feasibility ratio under this definition. 
The results show that for $\mathcal{P}_1$ and $\mathcal{P}_2$, the non-conditional approach offers a higher cumulative feasibility ratio than the conditional approach. 
As the problem complexity increases further (e.g., $\mathcal{P}_3$), the gap between the two approaches diminishes, resulting in comparable, and in some cases worse (e.g., $\Delta \tau =$ $1.0$ and $1.5$ days), cumulative feasibility ratios compared to the conditional approaches. 
However, the limited statistical sample size at this depth of robustness suggests potential uncertainty in these comparative metrics, necessitating careful interpretation of the observed trends.

\ifthenelse{\boolean{includefigures}}
{
    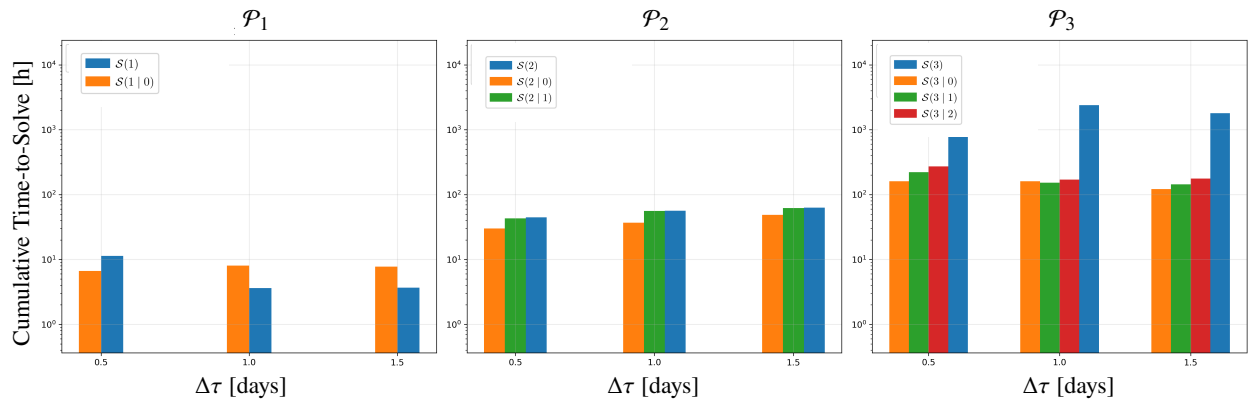
\begin{figure}[!htb]
        \centering
        \begin{tikzonimage}[keepaspectratio, width=0.99\textwidth]{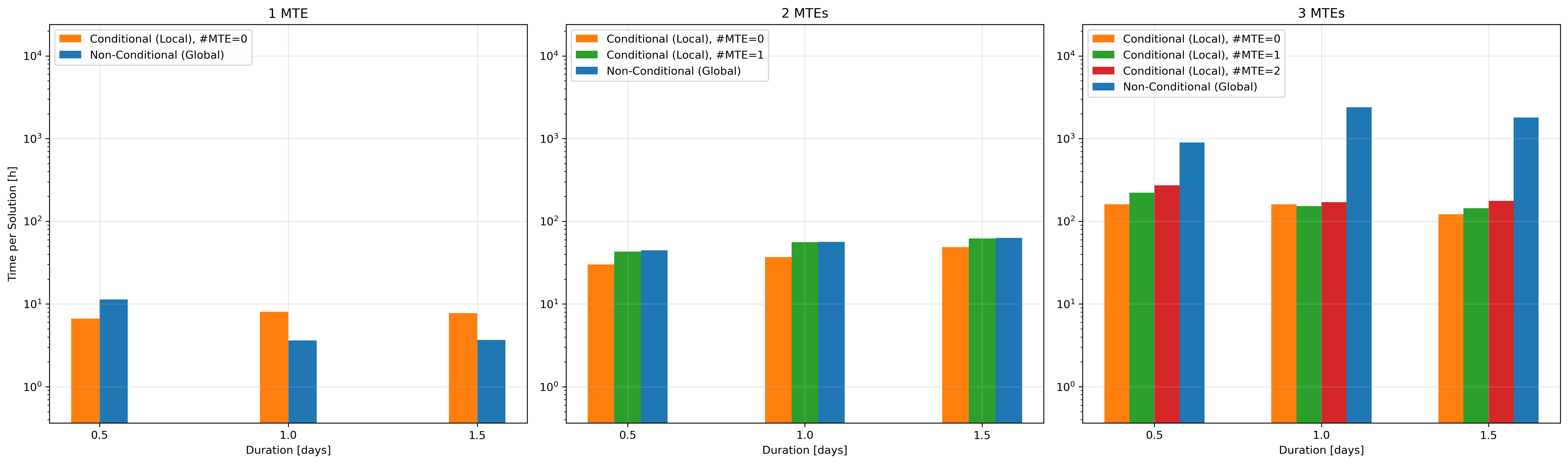}
            \small
            \node[fill=white, opacity=1.0, text opacity=1, anchor=south, rotate=90] at (0.015,0.50) {Cumulative Time-to-Solve [h]};
            \node[fill=white, opacity=1.0, text opacity=1, anchor=south] at (0.19,0.95) {$\mathcal{P}_1$};
            \node[fill=white, opacity=1.0, text opacity=1, anchor=south] at (0.515,0.95) {$\mathcal{P}_2$};
            \node[fill=white, opacity=1.0, text opacity=1, anchor=south] at (0.845,0.95) {$\mathcal{P}_3$};
            \node[fill=white, opacity=1.0, text opacity=1, anchor=south] at (0.1775,-0.060) {$\Delta\tau$ [days]};
            \node[fill=white, opacity=1.0, text opacity=1, anchor=south] at (0.5125,-0.060) {$\Delta\tau$ [days]};
            \node[fill=white, opacity=1.0, text opacity=1, anchor=south] at (0.8500,-0.060) {$\Delta\tau$ [days]};
            \node[fill=white,opacity=1.0,anchor=center] at (0.1025,0.855){\includegraphics[keepaspectratio, width=0.75 in]{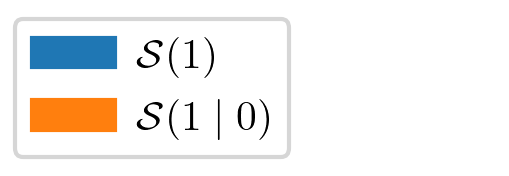}};
            \node[fill=white,opacity=1.0,anchor=center] at (0.43,0.83){\includegraphics[keepaspectratio, width=0.75 in]{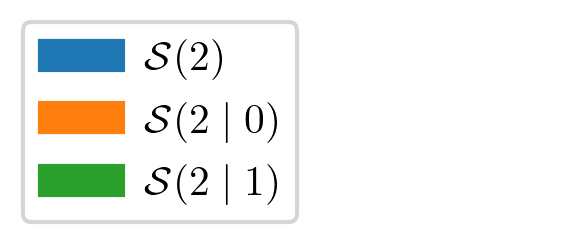}};
            \node[fill=white,opacity=1.0] at (0.76,0.81){\includegraphics[keepaspectratio, width=0.75 in]{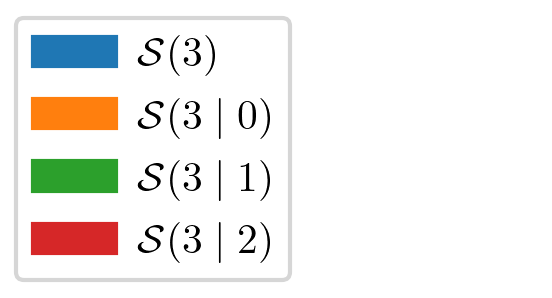}};
        \end{tikzonimage}
        \caption{Cumulative Average Time-to-Solve Comparison for Conditional and Non-Conditional Search Strategies}
        \label{fig: adjusted time to solve}
    \end{figure}
}
{
}

The solving time metric requires analogous adjustment for the conditional approaches.
Continuing with the $\mathcal{S}(1 | 0)$ example, the \emph{cumulative} solving time would be the sum of the $\mathcal{S}(0)$ solving time and the $\mathcal{S}(1 | 0)$ solving time. 
Figure \ref{fig: adjusted time to solve} presents these cumulative solving time metrics for the same strategies discussed earlier.
For $\mathcal{P}_1$, the non-conditional strategy out-performs the condition strategy. 
However, as the problem dimensionality increases, the conditional approach, specifically $\mathcal{S}(k | 0)$, exhibits superior performance, achieving reduced cumulative solving times relative to non-conditional search.
\section{Conclusion}
\label{section: conclusion}

This study introduced a novel initial guess generation strategy, the \emph{conditional global} approach, and compared it with a baseline na\"ive \emph{non-conditional global} approach, in the context of robust low-thrust design against missed-thrust events. 
A quantitative assessment using three key algorithmic performance metrics - feasibility ratios, average solving time, and solution quality - across varying robust problem complexity revealed distinct methodological characteristics.

The non-conditional approach, while enabling a more comprehensive exploration of the solution space, exhibited degraded performance at higher robustness depths, characterized by declining feasibility ratios and escalating average solving times. 
Conversely, the conditional approach, which leverages prior solutions to help narrow the search space, demonstrated higher feasibility ratios and lower average solving times, particularly when warm-started with previously solved non-robust solutions. 
Notably, this contradicts the intuition that warm-starting the search process with partially robust solutions would provide better initial guesses, suggesting that feasible region accessibility dominates partial robustness considerations in initial guess strategy selection.
When accounting for the complete computational framework, including the computational overhead in generating the solutions used as initial guesses for the conditional search, the relative performance advantages between methodologies demonstrate a dependence on the problem complexity. 
The non-conditional approach exhibits higher cumulative feasibility ratio at lower robustness depths, while conditional methods achieve better cumulative average solving times at higher robustness depths. 

These findings illuminate the fundamental trade-offs between exploration and exploitation in robust low-thrust trajectory design. 
While non-conditional search offers solution diversity without prior information, its effectiveness rapidly diminishes with increasing problem complexity. 
The conditional approach, through targeted refinement of existing solutions from simpler robust problems, maintains higher convergence rates without compromising solution quality. 
Future developments could enhance these methodologies through the investigation of alternate realization mapping strategies, exploration of problems with higher depths of robustness, and integration with advanced global optimization frameworks, potentially enabling more efficient robust mission design under missed thrust events.

\section*{Acknowledgments}
The simulations presented in this article were performed on computational resources managed and supported by Princeton Research Computing, a consortium of groups including the Princeton Institute for Computational Science and Engineering (PICSciE) and the Office of Information Technology's High Performance Computing Center and Visualization Laboratory at Princeton University.

\bibliography{references}


\end{document}